\documentclass{mcom-l}

\usepackage{graphicx}
\usepackage{amsmath}
\usepackage{amssymb}
\usepackage{amsfonts}
\usepackage{epsfig}
\usepackage{boxedminipage}
\usepackage{multirow}
\usepackage{booktabs}
\usepackage{graphicx}
\usepackage{color}
\usepackage{soul}

\newtheorem{theorem}{Theorem}[section]
\newtheorem{lemma}[theorem]{Lemma}
\newtheorem{corollary}[theorem]{Corollary}
\theoremstyle{definition}

\theoremstyle{remark}

\numberwithin{equation}{section}

\newcommand{\bflambda}{{\mbox{\boldmath $\lambda$}}}
\newcommand{\bfx}{{\mbox{\boldmath $x$}}}
\newcommand{\bfy}{{\mbox{\boldmath $y$}}}
\newcommand{\bft}{{\mbox{\boldmath $t$}}}
\newcommand{\bfw}{{\mbox{\boldmath $w$}}}
\newcommand{\bfu}{{\mbox{\boldmath $u$}}}
\newcommand{\bfv}{{\mbox{\boldmath $v$}}}

\begin{document}

\title[Solving Multiple-Block Separable Problems Using Two-Block ADMM]{Solving Multiple-Block Separable Convex Minimization Problems Using Two-Block Alternating Direction Method of Multipliers}


\author{Xiangfeng Wang}
\address{Department of Mathematics, Nanjing University, 22 Hankou Road, Nanjing, 210093, P.R. China.}
\email{xfwang.nju@gmail.com}

\author{Mingyi Hong}
\address{Department of Electrical and Computer Engineering, University of Minnesota, Minneapolis, 55455, USA.}
\email{mhong@umn.edu}

\author{Shiqian Ma}
\address{Department of Systems Engineering and Engineering Management, The Chinese University of Hong Kong, Shatin, N. T., Hong Kong.}
\email{sqma@se.cuhk.edu.hk}

\author{Zhi-Quan Luo}
\address{Department of Electrical and Computer Engineering, University of Minnesota, Minneapolis, 55455, USA.}
\email{luozq@ece.umn.edu}


\date{}

\dedicatory{}

\begin{abstract}
In this paper, we consider solving multiple-block separable convex minimization problems using alternating direction method of multipliers (ADMM). Motivated by the fact that the existing convergence theory for ADMM is mostly limited to the two-block case, we analyze in this paper, both  theoretically and numerically, a new strategy that first transforms a multi-block problem  into an equivalent two-block problem (either in the primal domain or in the dual domain) and then solves it using the standard two-block ADMM. In particular, we derive convergence results for this two-block ADMM approach to solve multi-block separable convex minimization problems, including an improved $\mathcal{O}(1/\epsilon)$ iteration complexity result. Moreover, we compare the numerical efficiency of this approach with  the standard multi-block ADMM on several separable convex minimization problems which include basis pursuit, robust principal component analysis and latent variable Gaussian graphical model selection. The numerical results show that the multiple-block ADMM, although lacks theoretical convergence guarantees, typically outperforms two-block ADMMs.
\end{abstract}

\maketitle


\bibliographystyle{amsplain}



\setcounter{equation}{0}
\section{Introduction}\label{secIntroduction}
In this paper, we consider the following convex optimization problem with $m$ block variables and the objective being the sum of $m\ (m\geq 2)$ separable convex functions:
\begin{eqnarray}\label{MultipleBlockProblem}
\min_{\{x_1,x_2,...,x_m\}}&& \sum_{i=1}^m f_i (x_i),\nonumber\\
\hbox{s.t.}&& \sum_{i=1}^m A_i x_i = b,\\
&& x_i \in {\mathcal{X}}_i\subseteq \mathbb{R}^{n_i},\ i = 1,\cdots,m,\nonumber
\end{eqnarray}
where $f_i:\mathbb{R}^{n_i}\rightarrow \mathbb{R}$ is a proper closed convex function (not necessarily smooth), ${\mathcal{X}}_i$ is a closed convex set and $A_i\in \mathbb{R}^{\ell\times n_i}$, $b\in \mathbb{R}^{\ell}$, and $\mbox{\boldmath$x$}=(x_1,\cdots, x_m)\!\in\!\prod_i \mathcal{X}_i$ is a partition of the decision variable.

The multi-block convex optimization  problem (\ref{MultipleBlockProblem}) arises naturally in many practical applications. For example, consider the Robust Principal Component Analysis (RPCA) problem \cite{CandesRPCA2011}, whereby the goal is to recover the low-rank and sparse components of a given matrix $M\in \mathbb{R}^{\ell\times n}$ by solving the following nonsmooth convex optimization problem
\begin{eqnarray}\label{RobustPCA}
\min_{L,S,Z} && \|L\|_* + \tau \|S\|_1,\\
\hbox{s.t.} && L + S + Z = M,\nonumber\\
&&Z\in {\mathcal{H}}:=\{Z\in \mathbb{R}^{\ell\times
n}\mid \|P_{\Omega}(Z)\|_F\le \delta\}.\nonumber
\end{eqnarray}
In the above formulation, $L \in \mathbb{R}^{\ell\times n}$ and $S \in \mathbb{R}^{\ell\times n}$ are respectively the low-rank and the sparse components of the matrix $M$, while $Z$ represents the observation noise. The notation $P_{\Omega}(\cdot)$ signifies the entry-wise projection operator for a given index set $\Omega$:
\[
P_{\Omega}(Z)=\left\{\begin{array}{ll}
Z_{ij},& \mbox{if } (i,j)\in \Omega,\\
0,&\mbox{otherwise,}
\end{array}\right.
\]
while $\|\cdot\|_*$, $\|\cdot\|_1$ and $\|\cdot\|_F$ denote respectively the matrix nuclear norm (i.e., the sum of the matrix singular eigenvalues), the $L_1$ and Frobenius norm of a matrix. Clearly problem (\ref{RobustPCA}) corresponds to the case of $m=3$ in problem \eqref{MultipleBlockProblem}, with $x_1 = L$, $x_2 = S$, $x_3 = Z$ and
$$
f_1(L) := \|L\|_*,\quad f_2(S) := \tau\|S\|_1,\quad f_3(Z) :=
{\mathcal{I}}_{\mathcal{H}}(Z),
$$
where $\mathcal{I}_{\mathcal{H}}(\cdot)$ denotes the indicator function for the set ${\mathcal{H}}$.

Similarly, the so-called latent variable Gaussian graphical model selection (LVGGMS) problem \cite{Chandrasekaran2012}, which  is closely related to the inverse covariance matrix estimation problem, is also in the form of \eqref{MultipleBlockProblem}. In particular, suppose $(X,Y)$ is a pair of $(p+r)$-dimensional joint multivariate Gaussian random variables, with covariance matrix denoted by $\Sigma_{(X,Y)} := [\Sigma_X,\Sigma_{XY};\Sigma_{YX}, \Sigma_Y ]$ and its inverse $\Theta_{(X,Y)} :=[ \Theta_X,\Theta_{XY};\Theta_{YX}, \Theta_Y ]$ respectively. The random variable $X:= (X_1,X_2,\cdots,X_p)^T$ is observable while $Y := (Y_1,Y_2,\cdots,Y_r)^T$ is the latent (or hidden) random variable. In many applications, we typically have $r\ll p$. Moreover, the marginal distribution of the observed variables $X$ usually follows a sparse graphical model and hence its concentration matrix $ \Theta_X$ is sparse. Notice that the inverse of the covariance matrix for $X$ can be expressed as
\begin{align}
\Sigma_X^{-1} = \Theta_X - \Theta_{XY} \Theta_Y^{-1} \Theta_{YX}\label{eqInverseCovariance},
\end{align}
which is the difference between the sparse term $\Theta_X$ and the low-rank term $\Theta_{XY} \Theta_Y^{-1} \Theta_{YX}$ (since $r$ is much less than $p$). Thus, the task of estimating the sparse marginal concentration matrix $\Theta_X$ can be accomplished by solving the following regularized maximum likelihood problem
\begin{eqnarray}\label{LVGMS}
\min_{S,L}&&\langle S - L, \hat{\Sigma}_X \rangle - \hbox{logdet}(S - L) + \alpha_1 \| S \|_1 + \alpha_2 {\mathbf{Tr}}(L),\\
\hbox{s.t.}&&S - L \succ 0,\ L\succeq 0,\nonumber
\end{eqnarray}
where $\hat{\Sigma}_X \in \mathbb{R}^{p\times p}$ is the sample covariance matrix of $X \in \mathbb{R}^{n\times p}$ and ${\mathbf{Tr}}(L)$ denotes the trace of matrix $L \in \mathbb{R}^{p\times p}$, while $S \in \mathbb{R}^{p\times p}$. Evidently, we can rewrite (\ref{LVGMS}) in the following equivalent form by introducing a new variable $R \in \mathbb{R}^{p\times p}$,
\begin{eqnarray}\label{LVGMSeuqal}
\min_{S,L}&&\langle R, \hat{\Sigma}_X \rangle - \hbox{logdet}(R) + \alpha_1 \| S \|_1 + \alpha_2 {\mathbf{Tr}}(L) + {\mathcal{I}}(L\succeq 0),\\
\hbox{s.t.}&&R = S - L,\nonumber
\end{eqnarray}
where the constraint $R\succ 0$ is implicitly imposed by having the term $-\hbox{logdet}(R)$ in the objective function. It is easily seen that LVGGMS corresponds to the three block case $m = 3$ in problem (\ref{MultipleBlockProblem}) with $x = (R, S, L)$ and
$$
f_1(R) := \langle R, \hat{\Sigma}_X \rangle - \hbox{logdet}(R),\  f_2(S) := \alpha_1 \|S\|_1,\  f_3(L) := \alpha_2 {\mathbf{Tr}}(L) + {\mathcal{I}}(L\succeq 0),
$$
where the linear constraint is $R - S + L = 0$.

Problem \eqref{MultipleBlockProblem} is a structured convex problem with a separable objective function and a single linear equality constraint. A popular algorithm for solving this class of problem is the so-called {\em{Alternating Direction Method of Multipliers}} (ADMM). To outline the basic steps of the ADMM, we first introduce the augmented Lagrangian function for problem \eqref{MultipleBlockProblem}
\begin{equation}\label{classicalLag}
L_{\beta} (x_1,\cdots,x_m,\lambda) = \sum\limits_{i=1}^m f_i (x_i) -
\langle \lambda, \sum\limits_{i=1}^m A_i x_i - b \rangle +
\frac{\beta}{2}\left\| \sum\limits_{i=1}^m A_i x_i - b \right\|^2,
\end{equation}
where $\beta$ is the penalty parameter for the violation of the linear constraint and $\langle \cdot,\cdot \rangle$ denotes the standard trace inner product. The ADMM method is a Gauss-Seidel iteration scheme in which the primal block variables $\{x_i\}$ and the Lagrangian multiplier $\lambda$ for the equality constraint are updated sequentially. Specifically, for a fixed stepsize $\beta> 0$, the ADMM for solving problem \eqref{MultipleBlockProblem} can be described as follows:
\begin{center}
\begin{tabular}{@{}llr@{}}\toprule
{\bf{\quad Algorithm 1: Alternating Direction Method of Multipliers for (\ref{MultipleBlockProblem})}}\\
\hline\\
\qquad Initialize $x_2^0,\cdots,x_m^0,\lambda^0$, and $\beta$.\\
\qquad For $k = 1,2,\cdots $, do\\
\qquad \qquad $\bullet$ Compute $x_i^{k+1}$, $\forall\  i = 1, \cdots,m$,\\
 \qquad\qquad$x_i^{k+1}=\arg\min\limits_{x_i\in {\mathcal{X}}_i} f_i(x_i) + \frac{\beta}{2}\left\| \sum\limits_{j=1}^{i-1}\! A_j x_j^{k+1} + A_i x_i + \sum\limits_{j=i+1}^m\!\! A_j x_j^k - b - \frac{\lambda^k}{\beta} \right\|^2,$\\
\qquad \qquad $\bullet$ Compute $\lambda^{k+1}$,\\
\qquad \qquad$\lambda^{k+1} = \lambda^k - \beta \left(\sum\limits_{i=1}^m\! A_i x_i^{k+1} - b\right).$\\[5pt]
\hline
\end{tabular}
\end{center}

The history of ADMM dates back to 1970s in \cite{GabayMercier1976,Glowinski1975} where the method was first developed for solving 2-block separable convex problems. In \cite{Gabay1983}, it is shown that ADMM can be interpreted as a special case of an operator splitting method called {\em{Douglas-Rachford Splitting Method}} (DRSM) for finding a zero of the sum of two monotone operators $\mathcal{A}$ and $\mathcal{B}$ \cite{DouglasRachford1956,LionsMercier1979}. Moreover, ADMM can also be related to Spingarn's method called {\em{Partial Inverse}} \cite{Spingarn1983,Spingarn1985}. Recently, ADMM has found its application in solving a wide range of large-scale problems from statistics, compressive sensing, image and video restoration, and machine learning, see e.g., \cite{BoydADMMsurvey2011,ChenHeYuan2012,GoldfarbMaScheinberg2012,HeXuYuan2011,MaXueZou2012, NgWeissYuan2010,XueMaZou2012,YangZhang2011,Yuan2012} and the references therein.

The ADMM convergence has long been established in the literature for the case of two block variables (i.e. $m=2$ in problem \eqref{MultipleBlockProblem}). References \cite{GabayMercier1976,Glowinski1975} show that the algorithm converges globally when each subproblem is solved exactly. Such convergence results have also been obtained in the context of DRSM. In \cite{EcksteinBertsekas1992}, the authors show that DRSM is a special case of the so-called {\em{Proximal Point Algorithm}} (PPA), for which the global convergence and the rate of convergence have been established (see \cite{Rockafellar1976}). Accordingly, under certain regularity conditions, the global convergence and the rate of convergence of DRSM (and hence ADMM) follow directly. On the other hand, for large scale problems such as those arising from compressive sensing, the global optimal solution for subproblems related to certain block variables may not be easily computable \cite{YangZhang2011}.  In these cases the classical ADMM needs to be modified accordingly so that those difficult subproblems are solved {\it inexactly}. In \cite{Eckstein1994,EcksteinBertsekas1992,HeLiaoHanYang2002,WangYuan2012,YangZhang2011,ZhangBurgerOsher2011}, the authors show that by performing a simple proximal gradient step for each subproblem, global convergence results similar to those for the classical ADMM can also be obtained. Beyond global convergence, there are a few results characterizing the iteration complexity and the convergence rate of the ADMM. The authors of \cite{HeYuan2012} have shown that to obtain an $\epsilon$-optimal solution, the worst-case iteration complexity of both exact and inexact ADMM is $\mathcal{O}(1/\epsilon)$, where the $\epsilon$ optimality is defined using both the constraint and objective violation. In \cite{Rockafellar1976}, Rockafellar has shown that if the inverse of the considered operator is Lipschitz continuous at origin point, the PPA converges linearly when the resolvent operator is solved either exactly or inexactly. Therefore the linear convergence of DRSM and ADMM follow directly under some assumptions on $\mathcal{A}$ and $\mathcal{B}$ in DRSM or $\{f_i\}$ and $\{A_i\}$ in ADMM. Further in \cite{LionsMercier1979}, Lions and Mercier have proved that when operator $\mathcal{B}$ is both coercive and Lipschitz, then DRSM converges linearly. Further in \cite{EcksteinBertsekas1990}, Eckstein and Bertsekas have shown the linear convergence rate of ADMM for linear programming. More recently the authors of \cite{DengYin2012} show that for both exact version and inexact version involving a proximal term, ADMM converges linearly if the objective function is strongly convex and its gradient is Lipschitz continuous in at least one block variable, and that the matrices $\{A_i\}$ satisfy certain rank assumptions.

However, when the number of block variables is greater than two, the convergence of ADMM has not been well understood. Recently, the authors of \cite{HongLuo2012} prove the global (linear) convergence of the ADMM for the multiple-block problem (\ref{MultipleBlockProblem}), under the following assumptions: {\it a)} for each $i$, $A_i$ is full column rank, {\it b)} $\beta$ is sufficiently small and {\it c)} certain error bounds hold for the problem \eqref{MultipleBlockProblem}. The full column rank condition a) can be dropped if each subproblem is solved inexactly. However when conditions {\it b)} and {\it c)} are not satisfied, even the global convergence of the algorithm is still open in the literature. As a result, a number of variants of the classical ADMM have been proposed; see \cite{HanYuanZhang2013,HeTaoYuan2012,HeTaoYuan2012further,Ma2012}. For example, in \cite{HeTaoYuan2012,HeTaoYuan2012further}, by adding an additional Gaussian back substitution correction step in each iteration after all the block variables are updated, the authors
establish the global convergence, iteration complexity and linear convergence rate (under some technical assumption on the iterates) of the modified algorithm. However such correction step is not always computationally efficient, especially when $A_i$'s are not efficiently invertible. In \cite{Ma2012}, an alternating proximal gradient method is proposed, in which the proximal version of the ADMM is applied to problem \eqref{MultipleBlockProblem}, after grouping the $m$ block variables into two blocks. However, the way that the block variables should be grouped is highly problem dependent. There is no general criterion guiding how such step should be done. Also recently, some multiple block splitting algorithms have been proposed for solving some models similar to \eqref{MultipleBlockProblem} in \cite{EcksteinSvaiter2009,GoldfarbMa2012}.

In this work, we systematically study ADMM algorithms for solving the  multi-block problem \eqref{MultipleBlockProblem}. We first propose two novel algorithms that apply the two-block ADMM to certain reformulation of the original multi-block problem. We show in detail how these algorithms are derived and analyze their convergence properties. We then report numerical results comparing the original multi-block ADMM with the proposed approaches on problems with multiple block structures such as the basis pursuit problem, robust PCA and latent variable Gaussian graphical model selection. Our numerical experiments show that the multi-block ADMM performs much better than the two-block ADMM algorithms as well as many other existing algorithms.

{\bf{Notation}}: Let $\| x \|_1 = \sum_{i=1}^n |x_i|$ and $\|x\|_2 = \sqrt{\sum_{i=1}^n x_i^2}$ denote the usual vector $\ell_1$-norm and $\ell_2$-norm respectively, where $x\in\mathbb{R}^n$ is a real vector. Let $\mathbf{I}_n\in\mathbb{R}^{n\times n}$ denote the $n\times n$ identity matrix. For a matrix $X$, let $\rho(X)$ denote its spectral radius. Throughout the paper we assume that the following {\em{proximity operator}} is easy to compute:
\begin{equation}\label{Proximal1}
\hbox{Prox}_{ g}^{\tau} (c): = \arg \min_{x\in {\mathcal{X}}}\ \tau g
(x) + \frac{1}{2}\| x - c \|^2,
\end{equation}
where $\tau>0$ and $c$ are given. For instance when $g(x) = \|x\|_1$,
\begin{equation}\label{Proximal2}
\left(\hbox{Prox}_{ \|\cdot\|_1}^{\tau} (c)\right)_i =
\left({\bf{shrinkage}}(c,\tau)\right)_i :=
\hbox{sgn}(c_i)\cdot\hbox{max}(|c_i| - \tau,0).
\end{equation}
We refer the readers to \cite{Combettes2011,Ma2012} for other easily computable proximity operators. The {\em{conjugate function}} of $f:\mathbb{R}^n\rightarrow \mathbb{R}$ is defined as
\begin{equation}\label{Conjugate}
f^* (y) := \sup_{x\in {\bf{dom}}f} (y^T x - f(x)).
\end{equation}
Let ${\mathcal{T}}$ be a set-valued operator, $c$ is any positive scalar, and denote the following operator
\begin{equation}\label{resolvent}
{\mathcal{J}}_{c{\mathcal{T}}}:= ({\mathcal{I}} + c {\mathcal{T}})^{-1},
\end{equation}
as the {\em{resolvent}} of ${\mathcal{T}}$.

The rest of this paper is organized as follows: the primal and dual splitting ADMM algorithms are presented and analyzed in Section \ref{subPrimalSplit} and Section \ref{subDualSplit}, respectively. In Section \ref{secNumerical}, numerical results for both synthetic and real problems are reported. Finally, some concluding remarks are given in Section \ref{secConclusion}.

\setcounter{equation}{0}
\section{The Proposed Algorithms}\label{secSplit}

Various modifications of the ADMM algorithm have been proposed in the literature
to deal with the multi-block problem \eqref{MultipleBlockProblem}.
Instead of proposing yet another ADMM variant, we propose to transform
any multi-block problem into an equivalent two-block
problem to which the classical two-block ADMM can be readily
applied. The main idea is to appropriately introduce some auxiliary
variables so that the original variables are completely decoupled
from each other. In this section, this technique will be explained
in details for both the primal and dual versions of problem
\eqref{MultipleBlockProblem}.

\subsection{Primal Splitting ADMM}\label{subPrimalSplit}

Introducing a set of auxiliary variables $\{y_i\}_{i=1}^{m}$ with $y_i=A_ix_i-\frac{b}{m}\in\mathbb{R}^{\ell}$ for all $i$, problem (\ref{MultipleBlockProblem}) can be expressed in the following equivalent form
\begin{eqnarray}\label{PrimalSplit}
\min_{\{x_i\}, \{y_i\}} && \sum_{i=1}^m f_i(x_i),\nonumber\\
\hbox{s.t.}&& A_i x_i - \frac{b}{m} = y_i, \ x_i\in {\mathcal{X}}_i,\ i = 1,\cdots,m,\\
&&\sum_{i=1}^m y_i = 0,\ y_i\in\mathbb{R}^{\ell}\nonumber.
\end{eqnarray}
Problems \eqref{MultipleBlockProblem} and \eqref{PrimalSplit} are equivalent in the sense that they share the same primal optimal solution set for the variables $\{x_i\}_{i=1}^{m}$, and achieve the same global optimal objective value. To apply the ADMM algorithm to the above reformulated problem, we write its (partial) augmented Lagrangian function, by keeping the constraint $\sum\limits_{i=1}^m y_i = 0$ and penalizing the rest of the constraints:
\begin{equation}\label{PrimalSplitLag}
L_{\beta}(\{x_i\},\{y_i\},\{\lambda_i\}) = \sum_{i=1}^m f_i(x_i) -
\sum_{i=1}^m \langle \lambda_i, A_i x_i - \frac{b}{m} - y_i \rangle
+ \frac{\beta}{2} \sum_{i=1}^m \left\| A_i x_i - \frac{b}{m} - y_i
\right\|^2,
\end{equation}
where $\lambda_i\in \mathbb{R}^{\ell}$ is the Lagrangian multiplier and $\beta$ is a penalty parameter. Denote $\mbox{\boldmath$\lambda$}: = (\lambda_1^T,\cdots,\lambda_m^T)^T(\in \Lambda:=\mathbb{R}^{\ell}\times\cdots\times\mathbb{R}^{\ell})$. Obviously, in \eqref{PrimalSplitLag}, all the $x_i$'s are separable with each other, so are all te $y_i$'s.  It is then natural to take $\mbox{\boldmath$x$}: = (x_1^T,\cdots,x_m^T)^T(\in\!\! {\mathcal{X}} := {\mathcal{X}}_1\times\cdots\times{\mathcal{X}}_m)$ and $\mbox{\boldmath$y$} : = (y_1^T \cdots,y_m^T)^T(\in {\mathcal{Y}} :=\{y\in\!\! \mathbb{R}^{\ell}\times\cdots\times\mathbb{R}^{\ell} \mid \sum_{i=1}^m y_i = 0\})$ as two block variables, and use the classical two-block ADMM to solve  (\ref{PrimalSplit}). The primal splitting ADMM is stated in the following table:
\begin{center}
\begin{tabular}{@{}llr@{}}\toprule
{\bf{\qquad\qquad Algorithm 2: Primal Splitting ADMM for (\ref{PrimalSplit})}}\\
\hline\\
\qquad Initialize $\{x_1^0,\cdots,x_m^0\},\{\lambda_1^0,\cdots,\lambda_m^0\}$, and $\beta$.\\
\qquad For $k = 0,1,2,\cdots $, do\\
\qquad \qquad $\bullet$ Compute $\{y_1^{k+1},\cdots,y_m^{k+1}\}$,\\
\qquad\qquad\qquad$\mbox{\boldmath$y$}^{k+1}=\arg\min\limits_{\mbox{\boldmath$y$}}\  \frac{\beta}{2}\sum\limits_{i=1}^m \left\| A_i x_i^k - \frac{b}{m} - y_i - \frac{\lambda_i^k}{\beta}\right\|^2,\ \hbox{s.t.}\ \sum\limits_{i=1}^m y_i = 0$,\\
\qquad \qquad $\bullet$ Compute $x_i^{k+1}$, $\forall i = 1, \cdots,m$,\\
\qquad\qquad\qquad$x_i^{k+1} = \arg\min\limits_{x_i\in {\mathcal{X}}_i}\  f_i (x_i) + \frac{\beta}{2}\left\| A_i x_i - \frac{b}{m} - y_i^{k+1} - \frac{\lambda_i^k}{\beta} \right\|^2,$\\
\qquad \qquad $\bullet$ Compute $\lambda_i^{k+1}$, $\forall i = 1, \cdots,m$,\\
\qquad\qquad\qquad$\lambda_i^{k+1} = \lambda_i^k - \beta \left(A_i x_i^{k+1} - \frac{b}{m} - y_i^{k+1}\right).$\\
\hline
\end{tabular}
\end{center}

We note that the subproblem for $\mbox{\boldmath$y$}$ is a
projection onto the hyperplane $\sum_iy_i=0$. As such, it admits the
following closed-form solution
\begin{equation}\label{PrimalSuby}
y_i^{k+1} = - \frac{1}{m}\left\{ \sum\limits_{i=1}^m A_i x_i^k - \frac{b}{m} - \frac{\lambda_i^k}{\beta} \right\} + \left( A_i x_i^k - \frac{b}{m} - \frac{\lambda_i^k}{\beta} \right),\ i=1,\cdots,m.
\end{equation}

Further, it is easy to see that the subproblem for $x_i$ can be solved efficiently by using (\ref{Proximal1}), provided that $A_i$ is an identity matrix (or any constant multiple of it). Else, $x_i$ can be updated by simply using a proximal gradient step:
\begin{eqnarray}\label{xLinearAlg2}
x_i^{k+1} \!\!\!\!&=&\!\!\!\! \arg\min\limits_{x_i\in {\mathcal{X}}_i} f_i (x_i) \!+\! \left\langle \beta A_i^T( A_i x_i^k \!-\! \frac{b}{m} \!-\! y_i^{k+1} \!-\! \frac{\lambda_i^k}{\beta}), x_i \!-\! x_i^k \right\rangle + \frac{\tau_i}{2}\left\| x_i \!-\! x_i^k \right\|^2,\nonumber\\
\!\!\!\!&=&\!\!\!\! \arg\min\limits_{x_i\in {\mathcal{X}}_i} f_i (x_i) + \frac{\tau_i}{2}\left\| x_i - x_i^k + \frac{\beta A_i^T( A_i x_i^k - \frac{b}{m} - y_i^{k+1}  - \frac{\lambda_i^k}{\beta})}{\tau_i} \right\|^2,\nonumber\\
\!\!\!\!&=&\!\!\!\! \hbox{Prox}_{f_i}^{\frac{1}{\tau_i}} \left(x_i^k - \frac{\beta A_i^T( A_i x_i^k - \frac{b}{m} - y_i^{k+1}  - \frac{\lambda_i^k}{\beta})}{\tau_i}\right),
\end{eqnarray}
where $\tau_i$ denotes the penalty parameter for the distance between $x_i^{k+1}$ and $x_i^k$.

Next we discuss the convergence of the above primal splitting ADMM algorithm.
\begin{theorem}\label{Algorithm1Convergence}
For any $\beta >0$, suppose the subproblem for $x_i$ is either exactly solved, or is solved inexactly using \eqref{xLinearAlg2} with $\tau_i > \beta \cdot \rho (A_i^T A_i)$. Let
$(\mbox{\boldmath$x$}^k, \mbox{\boldmath$y$}^k, \mbox{\boldmath$\lambda$}^k )$ be any sequence generated by Algorithm~2. Then starting with any initial point $(\mbox{\boldmath$x$}^0, \mbox{\boldmath$y$}^0, \mbox{\boldmath$\lambda$}^0 )\in {\mathcal{X}}\times{\mathcal{Y}}\times\Lambda$, we have
\begin{itemize}
\item[1.] The sequence $\{\mbox{\boldmath$\lambda$}^k\}$ converges to $\mbox{\boldmath$\lambda$}^*$, where $\mbox{\boldmath$\lambda$}^*$ is the dual optimal solution for problem \eqref{PrimalSplit};
\item[2.] The sequence $\{\sum\limits_{i=1}^m f_i (x_i^k)\}$ converges to $p^*$, where $p^*$ is the primal optimal value for problem \eqref{PrimalSplit};
\item[3.] The residual sequence $\{A_i x_i^k - y^k_i-\frac{b}{m}\}$ converges to 0 for each $i=1,\cdots,m$;
\item[4.] If the subproblem for $x_i$ is exactly solved for each $i$, then the sequence $\{A_i x_i^k\}$ and $\{y^k_i\}$ converge. Moreover, if $A_i$ has full column rank, then $\{ x_i^k \}$ converges.  When the subproblem for $x_i$ is solved inexactly using \eqref{xLinearAlg2} with $\tau_i > \beta \cdot \rho (A_i^T A_i)$, the sequence $\{(\bfx^k,\bfy^k)\}$ converges to  an optimal solution to problem
\eqref{PrimalSplit}.
\end{itemize}
\end{theorem}
{\bf{Proof}}. When the subproblems are solved exactly, we are actually using the classical two-block ADMM to solve the equivalent formulation (\ref{PrimalSplit}). As a result, the first three conclusions as well as the convergence of $\{ A_i x_i^k \}$ and $\{y^k_i\}$ follow directly from the classical analysis of the two-block ADMM (see, e.g., \cite[Section 3.2]{BoydADMMsurvey2011}). The convergence of $\{x^k_i\}$ is a straightforward consequence of the convergence of $\{ A_i x_i^k \}$ and the assumption that $A_i$'s are all full column rank. When the subproblems are solved inexactly via \eqref{xLinearAlg2}, because $\tau_i > \beta \cdot \rho (A_i^T A_i)$, all the conclusions are implied by the result in \cite[Theorem 1]{HeLiaoHanYang2002}. \hfill$\square$

Besides global convergence, we can elaborate on other convergence properties of {\bf Algorithm~2}. First, for both the exact case and the inexact proximal case with $\tau_i > \beta \cdot \rho (A_i^T A_i)$, we can obtain the following iteration complexity result by adopting the variational inequality framework developed in \cite[Theorem 4.1]{HeYuan2012}. To illustrate, let $\mbox{\boldmath$w$}:= (\mbox{\boldmath$y$}^T, \mbox{\boldmath$x$}^T,\mbox{\boldmath$\lambda$}^T )^T\in {\mathcal{Y}}\times {\mathcal{X}}\times\Lambda$, and let $\{\mbox{\boldmath$w$}^k\}$ denote the sequence generated by {\bf{Algorithm~2}}. Further we define
\begin{align}
\theta(\mbox{\boldmath$x$}) &:= \sum\limits_{i=1}^m f_i (x_i),\quad
F_1(\mbox{\boldmath$w$}) := \left(\begin{array}{c}
                            \mbox{\boldmath$\lambda$}\\
                            -\mbox{\boldmath$A$}^T\mbox{\boldmath$\lambda$}\\
                             A_1 x_1 - y_1 - \frac{b}{m}\\
        \vdots\\
        A_m x_m - y_m - \frac{b}{m}\\ \end{array}
                            \right),\nonumber
                            \end{align}
where $ \mbox{\boldmath$A$} := \hbox{diag}\left\{ A_1,\cdots,A_m \right\}$. By \cite[Theorem 4.1]{HeYuan2012}, we have that at any given iteration $K>0$,  the solution $\tilde{\mbox{\boldmath$w$}}^K=\frac{1}{K+1}\sum_{k=0}^{K}\bfw^k$ is an $\epsilon-$optimal solution for problem \eqref{MultipleBlockProblem}. That is, we have
\begin{equation}\label{Iterationcomplexity2}
\theta(\tilde{\mbox{\boldmath$x$}}^K) - \theta(\mbox{\boldmath$x$})
+ (\tilde{\mbox{\boldmath$w$}}^K - \mbox{\boldmath$w$})^T
F_1(\mbox{\boldmath$w$})\le \frac{C_p^1}{2(K+1)},\quad \forall\ \mbox{\boldmath$w$}\in
{\mathcal{Y}}\times{\mathcal{X}}\times\Lambda,
\end{equation}
where $C_p^1 = \max_{\bfw\in {\mathcal{W}}} \| \mbox{\boldmath$w$} - \mbox{\boldmath$w$}^0 \|_H^2$ and $H$ is a positive semi-definite matrix which is associated with $\{A_i\}$, $\beta$ and $\{\tau_i\}$. Note that at optimality, the left hand side of \eqref{Iterationcomplexity2} is no greater than zero, therefore the above inequality is indeed a possible measure of the
optimality gap, although it is implicit. In the following, we show explicitly that the objective values decrease at the rate ${\mathcal{O}}(1/K)$.

\begin{theorem}\label{ComplexityObjective}
Let $\{ \bfx^k, \bfy^k, \bflambda^k \}$ be the sequence generated by {\bf{Algorithm 2}}, and $(\bfx^*,\bfy^*,\bflambda^*)$ be any optimal solution, then we have
$$
\theta(\tilde{\bfx}^K) - \theta(\bfx^*) \le \frac{C_p^2}{2K},
$$
where $\tilde{\bfx}^K = \frac{1}{K}\sum\limits_{k=1}^K \bfx^k$ and $C_p^2 = \sum\limits_{i=1}^m \left\{ \frac{\|\lambda_i^0\|^2}{\beta} + \|x_i^* - x_i^0\|^2 + \beta \| y_i^* - A_i x_i^0 + \frac{b}{m} \|^2 \right\}$.
\end{theorem}
{\bf{Proof}}. Please see {\bf Appendix}.

It is worth noting that the complexity results presented in Theorem \ref{ComplexityObjective} and the one presented in \eqref{Iterationcomplexity2} do not imply each other. Moreover, in the proof of Theorem \ref{ComplexityObjective}, we have explored certain structure of {\bf Algorithm 2}, therefore this result does not carry over to the general ADMM algorithm.

Next we show that the linear rate of convergence for {\bf{Algorithm 2}} can also be established using existing results for the {\it Dauglas-Rachford Splitting Method} (DRSM).
To this end, we first derive the relationship between {\bf Algorithm 2} and the DRSM. {Recall that DRSM solves the following problem
\begin{equation}
\hbox{Find}\ \ u,\ \  \hbox{s.t.}\ \ 0\in {\mathcal{A}} (u) + {\mathcal{B}}(u), \label{eqDR}
\end{equation}
by generating two sequences $\{u^k\}$ and $\{v^k\}$  according to:
\begin{center}
\begin{tabular}{@{}llr@{}}\toprule
{\bf{\qquad\quad  Dauglas-Rachford Splitting Method}}\\
\hline\\
\qquad Initialize $\tau$, and $v^0, u^0 = J_{\tau {\mathcal{B}}} (v^0)$.\\
\qquad For $k = 1,2,\cdots $, do\\
\qquad \qquad $\bullet$ $v^{k+1} = J_{\tau {\mathcal{A}}} (2 J_{\tau {\mathcal{B}}} - {\mathcal{I}})(v^k) + ({\mathcal{I}} - J_{\tau {\mathcal{B}}})(v^k)$,\\
\qquad \qquad $\bullet$ $u^{k+1} = J_{\tau {\mathcal{B}}}(v^{k+1})$.\\
\\
\hline
\end{tabular}
\end{center}

To see the exact form of the operators ${\mathcal{A}}$ and ${\mathcal{B}}$ for {\bf Algorithm 2}, let us consider the dual formulation of \eqref{PrimalSplit}, stated below
\begin{eqnarray}\label{DRSMtransfer}
&& \max_{ \bflambda } \min_{ \bfx,\bfy \in {\mathcal{Y}} } \ \sum_{i=1}^m f_i(x_i) - \sum_{i=1}^m \langle \lambda_i, A_i x_i - \frac{b}{m} - y_i \rangle\nonumber\\
&=& \min_{\bflambda } \left\{ \max_{ \bfx } \ \sum_{i=1}^m \langle \lambda_i, A_i x_i \rangle -  f_i(x_i) + \max_{\bfy\in {\mathcal{Y}} } \sum_{i=1}^m \langle \lambda_i, - \frac{b}{m} - y_i \rangle\right\}\nonumber\\
&=& \min_{\bflambda } \left\{ {\mathcal{I}}(\bflambda: \lambda_1\!\! =\!\! \lambda_2 \!\!=\!\! \cdots \!\!=\!\! \lambda_m) + \sum_{i=1}^m \langle \lambda_i, - \frac{b}{m} \rangle + \sum_{i=1}^m f_i^* (A_i^T \lambda_i)\right\},
\end{eqnarray}
where ${\mathcal{I}}(\cdot)$ denotes the indicator function.  By setting \begin{align}
{\mathcal{A}}&: = \partial( {\mathcal{I}}(\bflambda: \lambda_1 \!\!=\!\! \lambda_2 \!\!=\!\! \cdots \!\!=\!\! \lambda_m) ) - (I,\cdots,I)^T\frac{b}{m}\label{eqA}\\
{\mathcal{B}}&: = \sum_{i=1}^m \partial (f_i^* \circ A_i^T)\label{eqB},
\end{align}
we can rewrite the dual form of \eqref{PrimalSplit} (i.e., eq. \eqref{DRSMtransfer}) equivalently as finding a $\bflambda^*$ that satisfies
\begin{align}
0\in {\mathcal{A}}(\bflambda^*) + {\mathcal{B}}(\bflambda^*) \label{eqDRADMM}.
\end{align}
Applying DRSM to solve \eqref{eqDRADMM}, we obtain the ($k+1$)th iterate as follows
\begin{eqnarray}
\!\!\!\!\!\!\bfv^{k+1} \!\!\!&=&\!\!\! J_{\beta {\mathcal{A}}} (2 J_{\beta {\mathcal{B}}} - {\mathcal{I}})(\bfv^k) + ({\mathcal{I}} - J_{\beta {\mathcal{B}}})(\bfv^k)\nonumber\\
\!\!\!&=&\!\!\!\arg\min_{\bfv}\!\left\{ {\mathcal{I}}(\bfv\!:\! v_i \!=\! v_j) \!-\! \sum_{i=1}^m \langle v_i, \frac{b}{m} \rangle \!+\! \frac{1}{2\beta}\| \bfv \!-\! (2 \bfu^k \!-\! \bfv^k) \|^2 \right\} \!+\! (\bfv^k \!-\! \bfu^k),\label{SubproblemV}
\end{eqnarray}
\begin{equation}
\bfu^{k+1} = J_{\beta {\mathcal{B}}}(\bfv^{k+1})=\arg\min_{\bfu}\ \sum_{i=1}^m f_i^* (A_i^T u_i) + \frac{1}{2\beta}\|\bfu - \bfv^{k+1}\|^2.\label{SubproblemU}
\end{equation}
Further, applying the {\em{Fenchel-Rockafellar Duality}} \cite[Definition 15.19]{BauschkeCombettes2011}, we can write the dual problem for \eqref{SubproblemV} and \eqref{SubproblemU} (with dual variables $\bfy$ and $\bfx$) as
\begin{eqnarray}
\bfy^{k+1} \!\!&=&\!\!\arg\min_{\bfy\in {\mathcal{Y}}}\left\{ \frac{\beta}{2} \sum_{i=1}^m \left\| - y_i \!-\! \frac{b}{m} \!-\! \frac{2 u_i^k \!-\! v_i^k}{\beta} \right\|^2 \right\}, \ v_i^{k+1} \!=\! u_i^k \!-\! \beta \left(- y_i^{k+1} \!-\! \frac{b}{m}\right), \nonumber\\
\bfx^{k+1} \!\!&=&\!\!\arg\min_{\bfx} \left\{ \sum_{i=1}^m f_i ( x_i) + \frac{\beta}{2}\left\| Ax_i - \frac{v_i^{k+1}}{\beta}\right\|^2 \right\},\ u_i^{k+1} = v_i^{k+1} - \beta A x_i^{k+1}\nonumber.
\end{eqnarray}
Then obviously when substituting $u_i^k = v_i^k - \beta A x_i^k$ into the subproblem about $\bfy$, we obtain
$$
\bfy^{k+1} = \arg\min_{\bfy\in {\mathcal{Y}}}\left\{ \frac{\beta}{2} \sum_{i=1}^m \left\| A x_i^k - y_i - \frac{b}{m}  - \frac{u_i^k}{\beta} \right\|^2 \right\}.
$$
Similarly, when substituting $v_i^{k+1} = u_i^k - \beta \left(- y_i^{k+1} - \frac{b}{m}\right)$ into the subproblem about $\bfx$, we can get the following equivalent problem for $\bfx$
$$
\bfx^{k+1} = \arg\min_{\bfx} \left\{ \sum_{i=1}^m f_i ( x_i) + \frac{\beta}{2}\left\| Ax_i - y_i^{k+1} - \frac{b}{m} - \frac{u_i^k}{\beta}\right\|^2 \right\}.
$$
Combining $u_i^{k+1} = v_i^{k+1} - \beta A x_i^{k+1}$ and $v_i^{k+1} = u_i^k - \beta \left(- y_i^{k+1} - \frac{b}{m}\right)$, we obtain the update of $\bfu$,
$$
u_i^{k+1} = u_i^k - \beta \left(A x_i^{k+1} - y_i^{k+1} - \frac{b}{m}\right).
$$
The above analysis indicates that the sequence $\{\bfu^k\}$ is the same as the multiplier sequence $\{\bflambda^k\}$ in {\bf{Algorithm 2}}. As a result, {\bf Algorithm 2} (or in general the two-block ADMM) can be considered as a special case of DRSM.

The linear convergence of DRSM has been well studied in \cite[{\sc{Proposition}} 4]{LionsMercier1979} with an assumption that operator ${\mathcal{B}}$ is both strongly monotone and Lipschitz, which means there exists $\alpha > 0$ and $M$ such that
$$
\|{\mathcal{B}} (x_1) - {\mathcal{B}} (x_2)\| \le M \| x_1 - x_2 \|,
$$
$$
\langle {\mathcal{B}} (x_1) - {\mathcal{B}} (x_2), x_1 - x_2 \rangle \ge \alpha \| x_1 - x_2 \|^2.
$$
By using the results in \cite{GoebelRockafellar2008,Rockafellar1996,RockafellarWets1998}, we can show that if $f_i$'s are all strongly convex with Lipschitz continuous gradients, and when $A_i$'s are all full row rank, then the operator ${\mathcal{B}}$ is strongly monotone and Lipschitz. As a result, the sequences $\{\bfx^k\}$, $\{\bfy^k\}$ and $\{\bflambda^k\}$ generated by {\bf{Algorithm 2}} converge linearly. Similarly, if each subproblem cannot be solved exactly, then the linear convergence of the inexact version of {\bf Algorithm 2} (cf. \eqref{xLinearAlg2}, with $\tau_i > \beta \cdot \rho (A_i^T A_i)$) can be established by following \cite[{\sc{Theorem}} 4]{DengYin2012}. Again we require that $f_i$'s are all strongly convex with Lipschitz continuous gradients, and $A_i$'s all have full row rank.

To this point, all the convergence results characterize the behavior of {\bf Algorithm 2} for solving problem \eqref{PrimalSplit}. As problem (\ref{MultipleBlockProblem}) is an equivalent reformulation of \eqref{PrimalSplit}, we can readily conclude that the sequence $\{\mbox{\boldmath$x$}^k\}$ generated by {\bf{Algorithm 2}} converges to the primal optimal solution of \eqref{MultipleBlockProblem}, if either each subproblem is exactly solved and $A_i$'s are all full rank, or the subproblems are solved using the proximal step \eqref{xLinearAlg2} with $\tau_i > \beta \cdot \rho (A_i^T A_i)$. Further, if $(\mbox{\boldmath$x$}^k, \bfy^k)$ is linearly convergent to an optimal solution of problem \eqref{PrimalSplit}, then $\mbox{\boldmath$x$}^k$ converges linearly to the primal optimal solution of (\ref{MultipleBlockProblem}).

\subsection{Dual Splitting ADMM}\label{subDualSplit} \setcounter{equation}{0}

We can also apply the splitting technique to the dual formulation (\ref{MultipleBlockProblem}) to derive a dual splitting ADMM algorithm. In particular, let us first write (\ref{MultipleBlockProblem}) in its saddle point form
\begin{equation}\label{saddlepointproblem}
\min_{x_i\in {\mathcal{X}}_i} \ \max_{\lambda}\  \sum_{i=1}^m f_i(x_i) - \langle \lambda, \sum_{i=1}^m A_i x_i - b \rangle.
\end{equation}
By exchanging the order of $\max$ and $\min$, and using the definition of the conjugate function of $f_i$, we can rewrite (\ref{saddlepointproblem}) equivalently as
\begin{eqnarray}\label{saddlepointproblem-rewrite}
&&\max_{\lambda}\ \min_{x_i\in {\mathcal{X}}_i}\  \sum_{i=1}^m f_i(x_i) - \langle \lambda, \sum_{i=1}^m A_i x_i - b \rangle\nonumber\\
&\Leftrightarrow&  \max_{\lambda}\ \min_{x_i\in {\mathcal{X}}_i} \left\{ \sum_{i=1}^m \left( f_i(x_i) - \langle A_i^T \lambda, x_i \rangle\right)\right\} + \langle \lambda, b \rangle\nonumber\\
&\Leftrightarrow& \max_{\lambda}\  - \sum_{i=1}^m f_i^* (A_i^T \lambda) + \langle \lambda, b \rangle\  \Leftrightarrow\ \min_{\lambda}\  \sum_{i=1}^m f_i^* (A_i^T
\lambda) - \langle \lambda, b \rangle,
\end{eqnarray}
where $\lambda$ denotes the dual variable of (\ref{MultipleBlockProblem}). We then split the dual variable $\lambda$ by introducing a set of auxiliary variables $\{\lambda_i\}_{i=1}^{m}$, and rewrite \eqref{saddlepointproblem-rewrite} as
\begin{eqnarray}\label{DualSplit}
\min_{\lambda,\ \lambda_i}&&\sum_{i=1}^m f_i^* (A_i^T \lambda_i) - \langle \lambda, b \rangle,\quad \hbox{s.t.}\  \lambda = \lambda_i,\ i = 1,\cdots,m.
\end{eqnarray}
It is obvious that each primal optimal solution $\{\lambda_i^*\}_{i=1}^m$, $\lambda^*$ of \eqref{DualSplit} corresponds to a dual optimal solution of \eqref{MultipleBlockProblem} ($\lambda_i^* = \lambda^*$). The augmented Lagrangian function for this dual problem can be expressed as follows:
\begin{equation}\label{DualLag}
L(\{\lambda_i\},\lambda,\{t_i\}) = \sum_{i=1}^m f_i^* (A_i^T \lambda_i) - \langle\lambda, b\rangle - \sum_{i=1}^m \langle t_i, \lambda - \lambda_i \rangle + \frac{\beta}{2}\sum_{i=1}^m\| \lambda - \lambda_i \|^2,
\end{equation}
where $\beta$ is the penalty parameter for the constraints violation, and $t_i\in \mathbb{R}^{\ell}$ is the Lagrangian multiplier associated with the constraint $\lambda = \lambda_i$. Denote $\mbox{\boldmath$t$}: = (t_1^T,\cdots,t_m^T )^T\in \mathbb{R}^{m\ell}$ and $\tilde{\mbox{\boldmath$\lambda$}}: = (\lambda_1^T,\cdots,\lambda_m^T )^T\in \mathbb{R}^{m\ell}$. It is clear now that optimizing the augmented Lagrangian for each auxiliary variable $\lambda_i$ is independent of all other auxiliary variables $\{\lambda_j\}_{j\ne i}$. Consequently by treating $\tilde{\mbox{\boldmath$\lambda$}}$ and $\lambda$ as two block variables, we can again apply the two-block ADMM to solve (\ref{DualSplit}). In the following we take a closer look at the structure of each subproblem.

The subproblem for $\tilde{\bflambda}$ is related to the conjugate function $f_i^*(\cdot)$. At the $k$th iteration, this subproblem can be explicitly expressed as the following $m$ independent problems (one for each variable $\lambda_i$):
\begin{equation}\label{SubLambdadual}
\lambda_i^{k+1} = \arg\min_{\lambda_i}\ f_i^*(A_i^T \lambda_i) +
\frac{\beta}{2}\left\| \lambda_i - \lambda^{k+1} + \frac{t_i^k}{\beta}
\right\|^2.
\end{equation}
By the classical {\em{Fenchel-Rockafellar}} duality \cite{BauschkeCombettes2011,Rockafellar1996},  the dual problem of (\ref{SubLambdadual}) can be expressed as
\begin{equation}\label{SubLambdaprimal}
x_i^{k+1} = \arg\min_{x_i\in {\mathcal{X}}_i}\ f_i(x_i) +
\frac{1}{2\beta}\left\| A_i x_i - \beta (\lambda^{k+1} -
\frac{t_i^k}{\beta}) \right\|^2,
\end{equation}
where $\{x_i\}$ is precisely the set of primal variables of (\ref{MultipleBlockProblem}). The relationship between $\lambda_i^{k+1}$ and $x_i^{k+1}$ is as follows
\begin{align}\label{eqLambda}
\lambda_i^{k+1} = \lambda^{k+1} - \frac{t_i^k}{\beta} - \frac{1}{\beta} A_i x_i^{k+1}.\end{align}
The dual splitting ADMM is stated formally in the following table.
\begin{center}
\begin{tabular}{@{}llr@{}}\toprule
{\bf{\qquad\quad Algorithm 3: Dual Splitting ADMM for (\ref{DualSplit})}}\\
\hline\\
\qquad Initialize $\{\lambda_1^0,\cdots,\lambda_m^0\},\{t_1^0,\cdots,t_m^0\}$, and $\beta$.\\
\qquad For $k = 0,1,2,\cdots $, do\\
\qquad \qquad $\bullet$ Compute $\lambda^{k+1}$,\\
\qquad\qquad\qquad$\lambda^{k+1}= \arg\min\limits_{\lambda}\  - \lambda^T b - \sum\limits_{i=1}^m \lambda^T t_i^k + \frac{\beta}{2}\sum\limits_{i=1}^m \| \lambda - \lambda_i^k \|^2$,\\
\qquad \qquad $\bullet$ Compute $\lambda_i^{k+1}$, for all $i = 1, \cdots,m$,\\
\qquad\qquad\qquad$x_i^{k+1} = \arg\min\limits_{x_i\in {\mathcal{X}}_i}\ f_i(x_i) + \frac{1}{2\beta}\left\| A_i x_i - \beta (\lambda^{k+1} - \frac{t_i^k}{\beta}) \right\|^2,$\\
\qquad\qquad\qquad $\lambda_i^{k+1} = \lambda^{k+1} - \frac{t_i^k}{\beta} - \frac{1}{\beta} A_i x_i^{k+1},$\\
\qquad \qquad $\bullet$ Compute $t_i^{k+1}$, $\forall \ i = 1, \cdots,m$,\\
\qquad\qquad\qquad$t_i^{k+1} = t_i^k - \beta (\lambda^{k+1} - \lambda_i^{k+1})$.\\
\hline
\end{tabular}
\end{center}

Similar to the case of primal splitting, the subproblem of $\lambda$ can be solved easily in closed-form
\begin{equation}\label{DualSuby}
\lambda^{k+1} = \frac{1}{m\beta}\left( b+ \sum_{i=1}^m (t_i^k + \beta \lambda_i^k) \right).
\end{equation}
Furthermore, the subproblem for the block variable $x_i$ can be solved efficiently if $A_i$ is an identity matrix (or any of its constant multiples), because (\ref{SubLambdaprimal}) can be efficiently solved by computing the proximity operator (\ref{Proximal1}). Else, a proximal gradient step can be performed, i.e.,
\begin{eqnarray}\label{xLinearAlg3}
x_i^{k+1} &=&\!\!\! \arg\min_{x_i\in {\mathcal{X}}_i}\ f_i(x_i) + \left\langle \frac{1}{\beta}A_i^T( A_i x_i^k - \beta \lambda^{k+1} + t_i^k),x_i - x_i^k \right\rangle + \frac{\tau_i}{2} \|x_i - x_i^k\|^2,\nonumber\\
&=&\!\!\! \arg\min_{x_i\in {\mathcal{X}}_i}\ f_i(x_i) + \frac{\tau_i}{2} \left\|x_i - x_i^k + \frac{1}{\tau_i\beta}A_i^T( A_i x_i^k - \beta \lambda^{k+1} + t_i^k)\right\|^2,\nonumber\\
&=&\!\!\! \hbox{Prox}_{f_i}^{\frac{1}{\tau_i}} \left(x_i^k - \frac{1}{\tau_i\beta}A_i^T( A_i x_i^k - \beta \lambda^{k+1} + t_i^k)\right),
\end{eqnarray}
where $\tau_i$ denotes the penalty parameter for the distance between $x_i^{k+1}$ and $x_i^k$.

Again, the global convergence of {\bf{Algorithm 3}} is a straightforward consequence of the standard convergence results for the
two-block ADMM.
\begin{theorem}\label{Algorithm2Convergence}
For any $\beta >0$, suppose the subproblem for $x_i$ is either exactly solved, or is solved using \eqref{xLinearAlg3} with $\tau_i > \frac{\rho (A_i^T A_i)}{\beta}$ associated with $\lambda_i$. Let $(\tilde{\mbox{\boldmath$\lambda$}}^k, \lambda^k, \mbox{\boldmath$t$}^k )$ be any sequence generated by Algorithm 3. Then starting with any initial point $(\tilde{\mbox{\boldmath$\lambda$}}^0, \lambda^0, \mbox{\boldmath$t$}^0 )\in \mathbb{R}^{m\ell}\times \mathbb{R}^{\ell}\times \mathbb{R}^{m\ell}$, we have
\begin{itemize}
\item[1.] The sequence $\{\mbox{\boldmath$t$}^k\}$ converges to  the dual optimal solution for
problem \eqref{DualSplit}.
\item[2.] The sequence $\left\{\sum\limits_{i=1}^m f_i^* (A_i^T \lambda_i^k) - \langle \lambda^k, b \rangle\right\}$ converges to the primal optimal value for problem
\eqref{DualSplit}.
\item[3.] The residual sequence $\{\lambda^k - \lambda_i^k\}$ converges to 0 for each $i=1,\cdots,m$.
\item[4.] The sequence $\{\tilde{\mbox{\boldmath$\lambda$}}^k, \lambda^k\}$ converges to the primal optimal solution for problem \eqref{DualSplit}.
\item[5.] For each $i=1,\cdots,m$, if the subproblem about $x_i$ is exactly solved, then the sequence $\{A_i x_i^k\}$ converges. If $A_i$ has full column rank, then $\{ x_i^k \}$ converges to $x_i^*$, for all $i=1,\cdots,m$; the same is true if the subproblem about $x_i$ is solved using \eqref{xLinearAlg3} with $\tau_i > \frac{\rho (A_i^T A_i)}{\beta}$.
\end{itemize}
\end{theorem}
{\bf{Proof}}. When the subproblems are solved exactly, {\bf Algorithm 3} corresponds to the classical two-block ADMM applied to solve the equivalent formulation (\ref{DualSplit}). As a result, the first four conclusions follow directly from the classical analysis of the two-block ADMM (see, e.g., \cite[Section 3.2]{BoydADMMsurvey2011}). In the last conclusion, the convergence of $\{ A_i x_i^k \}$ follows from the convergence of $\{\lambda^k_i\}$, $\lambda^t$ and $\{t^t_i\}$; see \eqref{eqLambda}. When the subproblems are solved inexactly  with $\tau_i > \frac{\rho (A_i^T A_i)}{\beta}$, there is no existing result which covers the convergence of the algorithm. We will provide a  proof for this case in the {\bf{Appendix}}. \hfill$\square$

Let us discuss some additional convergence properties of {\bf Algorithm 3}. First of all, it is possible to derive the iteration complexity for both the exact and the inexact versions of the dual splitting ADMM algorithm. For the exact version, its iteration complexity based on variational inequalities follows from the existing results \cite{HeYuan2012}. The iteration complexity of the inexact dual splitting ADMM algorithm is not covered by any existing result. As a result, in {\bf{Appendix}}, we provide a unified iteration complexity analysis for the dual splitting ADMM algorithm. Additionally, similar to Theorem \ref{ComplexityObjective}, we have the following result that bounds the gap of objective value. The proof is similar to that of Theorem \ref{ComplexityObjective}, thus we omit it for brevity.
\begin{corollary}\label{corComplexity}
Let $\{\bfx^k,\lambda^k\}$ be the sequence generated by {\bf{Algorithm 3}} (using either the exact minimization or the inexact version with $\tau_i > \frac{\rho (A_i^T A_i)}{\beta}$) and $\bfx^*$ be any primal optimal solution. Define $\tilde{\bfx}^K := \frac{1}{K} \sum_{k=1}^K \bfx^k$, $\tilde{\lambda}^K := \frac{1}{K} \sum_{k=1}^K \lambda^k$. We have
\begin{equation}\label{Iterationcomplexity3}
\theta(\tilde{\mbox{\boldmath$x$}}^K) - \theta(\mbox{\boldmath$x$})
+ \!\left(\!\begin{array}{c}\tilde{\bfx}^K - \mbox{\boldmath$x$}\\ \tilde{\lambda}^K - \lambda \end{array}\!\right)\!^T
F_2(\bfx,\lambda)\le \frac{C_d^1}{2K},\ \forall \!\left(\!\!\begin{array}{c} \mbox{\boldmath$x$}\\  \lambda \end{array}\!\!\right)\!\!\in {\mathcal{X}}\times{\mathbb{R}}^{\ell},
\end{equation}
$$
\theta(\tilde{\bfx}^K) - \theta(\bfx^*)\le \frac{C_d^2}{2K},
$$
$$
C_d^1 \!=\! \max\limits_{ \bfx\in {\mathcal{X}},\lambda\in {\mathbb{R}}^l} \left\{ \sum\limits_{i=1}^m \| x_i \!-\! x_i^0 \|_{G_i}^2 \!+\! \beta \| \lambda \!-\! \lambda_i^0 \|^2 \right\},
$$
$$
C_d^2 \!=\! \left\{ \frac{\beta}{2} \sum_{i=1}^m \left\| \lambda_i^0 \right\|^2 \!+\! \frac{1}{2\beta}\sum_{i=1}^m \left\| t_i^* \!-\! t_i^0 \right\|^2 \!+\! \frac{1}{2}\sum_{i=1}^m \left\| x_i^* \!-\! x_i^0 \right\|_{\tilde{P}_i}^2 \right\},
$$
where $F_2 (\bfx,\lambda) \!\!=\!\! {\small \!\left(\!\begin{array}{c} -A_1^T \lambda\\ \cdots \\ -A_m^T \lambda\\ \sum_{i=1}^m A_i x_i - b \end{array}\!\right)\!}$ and $G_i \!=\! \frac{1}{\beta} A_i^T A_i$, $\tilde{P}_i = 0$ for the exact version and $G_i \!=\! \tau_i I$, $\tilde{P}_i = \tau_i I - \frac{1}{\beta}A_i^T A_i$ for the inexact version.
\end{corollary}

Further, the linear rate of convergence for {\bf{Algorithm 3}} can also be established using existing results. We apply DSRM with two operators ${\mathcal{A}}\!:=\!\partial ({\mathcal{I}}({\sum_{i=1}^m t_i \!+\! b \!=\! 0}))$ and ${\mathcal{B}}:= \sum_{i=1}^m \partial (f_i \circ A_i^T)$. By following an analysis similar to that of {\bf{Algorithm 2}} and using assumptions that $f_i$'s are strongly convex and have Lipschitz continuous gradients, and that the matrices $A_i$'s have full row rank, we can prove that $\{\tilde{\bflambda}^k\}$, $\{\lambda^k\}$, $\{\bft^k\}$ and $\{\bfx^k\}$ converge linearly.

From the equivalence relationship between the problems \eqref{MultipleBlockProblem} and \eqref{DualSplit}, we can readily claim that the primal optimal solution $\lambda^*$ of (\ref{DualSplit}) is the dual optimal solution of (\ref{MultipleBlockProblem}).  Recall that from the discussion following \eqref{SubLambdaprimal}, $\{x_i\}$ is the set of primal variables for the original problem \eqref{MultipleBlockProblem}.

\subsection{Discussions}
Several existing methods for solving \eqref{MultipleBlockProblem} are similar to the two algorithms ({\bf{Algorithm 2}} and {\bf{Algorithm 3}}) proposed in this paper. In particular, Spingarn \cite{Spingarn1985} applied a method called {\em{partial inverse method}} \cite{Spingarn1983} to separable convex problems. This method can be directly applied to \eqref{MultipleBlockProblem} as follows. Let us define two subspaces $A$ and $B$ as:
$$
A :=\left \{ (\bfx,\bfu)\mid \sum_{i=1}^m u_i = 0\right\},\qquad B :=\left \{(0,\bfu)\mid u_1 = \cdots = u_m\right\},
$$
where $\bfx\in \mathbb{R}^n$, $u_i\in \mathbb{R}^{\ell}$, $\bfu:=(u_1^T,\cdots,u_m^T)^T\in \mathbb{R}^{m\ell}$.
Define the function
$$
F(\bfx,\bfu) = \left\{\begin{array}{l}
            \sum\limits_{i=1}^m f_i(x_i),\quad \hbox{if } A_i x_i - \frac{b}{m} = u_i, x_i\in {\mathcal{X}}_i, i = 1,\cdots,m,\nonumber\\
            +\infty,\qquad\quad \hbox{otherwise}.\nonumber\\
            \end{array}\right.
$$
Then problem \eqref{MultipleBlockProblem} is equivalent to the one that minimizes $F(\bfx,\bfu)$ over $A$. Define two operators $$
P_{A} (\bfx,\bfu) = \left(\bfx, \bfu - (I,\cdots,I)^T\sum\limits_{i=1}^m \frac{u_i}{m}\right),\quad P_{B} (\bfx,\bfu) = \left(0,(I,\cdots,I)^T \sum\limits_{i=1}^m \frac{u_i}{m}\right).
$$
To solve \eqref{MultipleBlockProblem}, partial inverse method generates a sequence of iterates $\{(\bfx^k,\bfy^k)\}$ and $\{(0,\bflambda^k)\}$:
$$
(\bfx^{k+1},\bfy^{k+1}) = P_{A}(\tilde{\bfx}^k,\tilde{\bfy}^k),\quad (0,\bflambda^{k+1}) = P_{B}(0,\tilde{\bflambda}^k),
$$
where $\{(\tilde{\bfx}^k,\tilde{\bfy}^k),(0,\tilde{\bflambda}^k)\}$ satisfies
\begin{eqnarray}
(\bfx_k,\bfy^k) + (0,\bflambda_k) &=& (\tilde{\bfx}^k,\tilde{\bfy}^k) + (0,\tilde{\bflambda}^k),\nonumber\\
\frac{1}{c_k} P_{A}(0,\tilde{\bflambda}^k) + P_{B}(0,\tilde{\bflambda}^k)&\in& \partial F \left( P_{A}(\tilde{\bfx}^k,\tilde{\bfy}^k) + \frac{1}{c_k}P_{B}(\tilde{\bfx}^k,\tilde{\bfy}^k) \right).\nonumber
\end{eqnarray}
with positive sequence $\{c_k\}$ bounded away from zero.

From \cite[Algorithm 2]{Spingarn1985}, it is known that the partial inverse method is the same as ADMM for minimizing $F(\bfx,\bfy)$ over $A$. That is, the variables $(\bfx,\bfy)$ in the {{Partial Inverse}} are the same as those in {\bf{Algorithm 2}} when the subproblems about $x_i$ are solved exactly.} However, the subproblem about $\bflambda$ in the partial inverse method additionally requires every component of $\bflambda$ to be equal, which is different from that in {\bf{Algorithm~2}}. Notice that the results in \cite{Spingarn1985} do not apply to the case when the subproblems for $x_i$ are solved inexactly.

{\bf Algorithm 3} is related to the proximal decomposition method proposed in \cite{CombettesPesquet2008}. The latter solves
\begin{align}\label{ProximalDecomposition}
\min_{x\in {\mathcal{X}}}& \quad \sum_{i=1}^m f_i (x)
\end{align}
by applying the  two block ADMM to the following reformulation
\begin{align}\label{ProximalDecomposition}
\min_{\{ x_1,x_2,\cdots,x_m \}}& \quad \sum_{i=1}^m f_i (x_i)\\
\mbox{s.t.}&\quad  x_i=y ,\  x_i\in {\mathcal{X}}, \ i=1,\cdots,m,
\end{align}
where ${\mathcal{X}}$ is the common closed convex constrained set for all $\{x_i\}_{i=1}^{m}$. In this decomposition, a single variable $x$ is split into $m$ copies $\{x_i\}_{i=1}^{m}$, and the consistency among these copies are enforced using the linking variable $y$. This decomposition technique is also used in {\bf Algorithm 3}, but for solving the dual reformulation of \eqref{MultipleBlockProblem}.

{\bf Algorithm 2} is closely related to the {\it distributed sharing} algorithm presented in \cite[Chapter 7]{BoydADMMsurvey2011}. Consider the following sharing problem, in which $m$ agents jointly solve the following problem
\begin{eqnarray}\label{Sharing}
\min_{\{ x_1,x_2,\cdots,x_m \}}&&\sum_{i=1}^m f_i (x_i) + g\left(\sum_{i=1}^m A_i x_i - b \right),\\
\hbox{s.t.}&&x_i \in {\mathcal{X}}_i,\ i = 1,\cdots,m,\nonumber
\end{eqnarray}
where $f_i(\cdot)$ is the cost related to agent $i$; $g(\cdot)$ is the cost shared among all the agents. The distributed sharing algorithm introduces a set of extra variables $y_i=A_i x_i - \frac{b}{m},\ \forall \ i$, and applies the two-block ADMM to the following reformulation
\begin{eqnarray}
\min_{\{ x_1,x_2,\cdots,x_m \}}&&\sum_{i=1}^m f_i (x_i) + g\left(\sum_{i=1}^m y_i \right),\nonumber\\
\hbox{s.t.}&& A_i x_i - \frac{b}{m} = y_i, \ x_i\in {\mathcal{X}}_i, \ i = 1,\cdots,m\nonumber.
\end{eqnarray}
To see the relationship between {\bf Algorithm 2} and the distributed sharing algorithm, we note that problem \eqref{MultipleBlockProblem}
is a special case of problem \eqref{Sharing}, with $g(\cdot)$ being the indicator function. Hence {\bf Algorithm 2} with the subproblems
being solved exactly can be viewed as a special case of the distributed sharing algorithm.

\section{Numerical Experiments}\label{secNumerical}
\setcounter{equation}{0}

In this section, we test {\bf Algorithms 1, 2 and 3}  on three problems: Basis Pursuit, Latent Variable Gaussian Graphical Model Selection and Robust Principal Component Analysis, and compare their performance with the {\em{Alternating Direction Method}}  (ADM)\cite{YangZhang2011}, {\em{Proximal Gradient based ADM}} (PGADM) \cite{MaXueZou2012}, {\em{ADMM with Gaussian Back Substitution}} (ADMGBS) \cite{HeTaoYuan2012}, and {\em{Variant Alternating Splitting Augmented Lagrangian Mehtod}}  (VASALM) \cite{TaoYuan2011}. Our codes were written in {\sc{Matlab 7.14}}(R2012a) and all experiments were conducted on a laptop with Intel Core 2 Duo@2.40GHz CPU and 4GB of memory.

\subsection{Basis Pursuit}
Consider the following basis pursuit (BP) problem \cite{ChenDonohoSaunders1998}
\begin{equation}\label{BP}
\min\  \|x\|_1\quad \hbox{s.t.}\quad A x = b,
\end{equation}
where $x\in \mathbb{R}^p$, and $A\in \mathbb{R}^{n\times p}$, $b\in \mathbb{R}^n$. This model has applications in compressed sensing where a sparse signal $x$ needs to be recovered using a small number of observations $b$ (i.e., $n\ll p$) \cite{ChenDonohoSaunders1998}.

By letting $\bfx = (\bfx_1,\cdots,\bfx_m)$, the BP problem can be viewed as a special case of problem \eqref{MultipleBlockProblem} with $m$ block variables. If we set $m=p$ (i.e., each component $x_i$ is
viewed as a single block variable), then {\bf Algorithm 1} can be used, with each of its primal iteration given by
\begin{eqnarray}
x_i^{k+1}&=&\arg\min_{x_i}\ \left|x_i\right| + \frac{\beta}{2}\left\|\sum_{j=1}^{i-1} a_j x_j^{k+1} + a_i x_i + \sum_{j = i+1}^p a_j x_j^k -b - \frac{\lambda^k}{\beta}\right\|^2\nonumber\\
&=& {\bf{shrinkage}}\left( \frac{1}{\sqrt{a_i^T a_i}} \left[ b + \frac{\lambda^k}{\beta} - \sum_{j=1}^{i-1} a_j x_j^{k+1} - \sum_{j = i+1}^p a_j x_j^k \right], \frac{1}{\beta (a_i^T a_i)} \right).
\end{eqnarray}
Alternatively, when the number of blocks $m$ is chosen as $m<p$, the primal subproblem in {\bf Algorithm 1} cannot be solved exactly. In this case, the inexact version of {\bf{Algorithm 2}} and {\bf{Algorithm 3}} can be used, where the primal subproblems \eqref{xLinearAlg2} and \eqref{xLinearAlg3} are respectively given by
\begin{eqnarray}
\eqref{xLinearAlg2}&\Leftrightarrow&\bfx_i^{k+1} = {\bf{shrinkage}}\left(\bfx_i^k - \frac{\beta A_i^T( A_i \bfx_i^k - \frac{b}{m} - \bfy_i^{k+1} - \frac{\lambda_i^k}{\beta})}{\tau_i},\frac{1}{\tau_i}\right),\nonumber\\
\eqref{xLinearAlg3}&\Leftrightarrow&\bfx_i^{k+1} = {\bf{shrinkage}}\left(\bfx_i^k - \frac{1}{\tau_i\beta}A_i^T( A_i \bfx_i^k - \beta \bflambda^{k+1} +\bft_i^k),\frac{1}{\tau_i}\right).\nonumber
\end{eqnarray}

In the following, we compare {\bf{Algorithm 1}} (with $m = p$), and the inexact versions of {\bf{Algorithm 2}}, {\bf{Algorithm 3}} (with $m = 2,5,10,20,50,100,200$) with the ADM algorithm \cite{YangZhang2011}, which has been shown to be effective for solving BP. {\bf Algorithms 1--3} are denoted as MULTADMM, PSADMM and DSADMM, respectively.

In our experiment, the matrix $A$ is randomly generated using standard Gaussian
distribution per element; the true solution $x^*$ is also generated
using standard Gaussian distribution, with $6\%$ sparsity level,
i.e., $94\%$ of the components are zero; the ``observation" vector $b$
is computed by $b = Ax^*$. The stepsize
$\beta$ is set to be $\frac{400}{\|b\|_1}$, $\frac{400}{\|b\|_1}$,
10 and $\frac{400}{\|b\|_1}$ for MULTADMM, PSADMM, DSADMM and ADM
respectively. Note that the stepsize for the DSADMM is chosen differently because it is the ADMM applied to the dual of \eqref{MultipleBlockProblem}, while the rest of the algorithms applies directly to the primal version of \eqref{MultipleBlockProblem}. To ensure convergence, the proximal parameters are set to be
$\tau_i =
1.01\beta\times\rho(A_i^T A_i)$ for PSADMM, $\tau_i =
1.01\times\frac{\rho(A_i^T A_i)}{\beta}$ for DSADMM, and $\tau =
1.01\beta\times\rho(A^T A)$ for ADM, respectively.
\begin{figure}[h]
\centering
\includegraphics[scale=0.4]{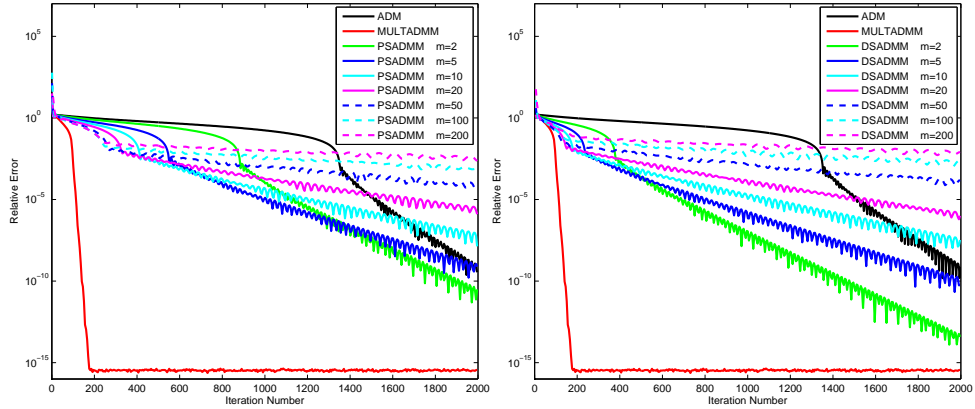}
\caption{Numerical performance of PSADMM and DSADMM for $n = 300$ and $p = 1000$}
\end{figure}

Figure 1 shows the convergence progress of the different algorithms. The curves in the figure represent the relative error for different algorithms along the iterations averaged over 100 runs. For a given iterate $x^k$, the
relative error is defined as
\begin{equation}\label{RelativeError}
\hbox{Relative Error} := \frac{\|x^k - x^*\|_2}{\|x^*\|_2}.
\end{equation}
The left part of Figure 1 shows the performance of MULTADMM, PSADMM and ADM. We observe that PSADMM converges faster than ADM when the number of blocks is relatively small. When the number of blocks increases, PSADMM converges fast at the beginning, but becomes slower after about $200$ iterations. This is because larger number of blocks results in smaller proximal parameter $\tau_i$, hence larger stepsizes can be taken \footnote{Note that the 2-norm of a matrix  (i.e., its largest singular value) cannot decrease if some of its columns are removed.}. On the other hand, it becomes increasingly difficult to simultaneously satisfy all the constraints. (Note that the number of constraint is the same as the number of block variables). We also observe that the MULTADMM performs much better than all other methods, in terms of both the convergence speed and the solution accuracy. The main reason for its superior performance is the exact solvability of the primal subproblems. Similar observations can be obtained from the right part of Figure 1, where the performance of MULTADMM, ADM and DSADMM are compared.

In Table 1 and Table 2, we report the performance of different algorithms for cases with $n=300,\ p=1000$ and $n=600,\ p=2000$, respectively. In these tables, `Iter' denotes the iteration number with $2000$ as the default maximum iteration number; `Obj' denotes the objective value; `Time' denotes the CPU time used; `$\sim$' indicates that the algorithm did not converge in 2000 iterations. Again we observe that to obtain the same accuracy, MULTADMM  requires significantly fewer iterations and less CPU time than all other methods. In the meantime, PSADMM and DSADMM perform better than ADM when the number of blocks is not too large.

\begin{center}
{\scriptsize\vskip 0.1cm
\begin{tabular}{|l|| c| c| c| c|| c| c| c| c|}
\multicolumn{9}{c}{\small{Table 1.\ Numerical comparison for Basis Pursuit with $n=300$ and $p=1000$}}\\
\hline
   \multicolumn{1}{|c||}{}&\multicolumn{4}{c||}{$\hbox{tol}=$1e-3}&\multicolumn{4}{c|}{$\hbox{tol}=$1e-5} \\
\hline
   \multicolumn{1}{|c||}{\hbox{Method}}&\hbox{Iter}&\hbox{Obj}&\hbox{Time}&\hbox{Error}&\hbox{Iter}&\hbox{Obj}&\hbox{Time}&\hbox{Error} \\
\hline
ADM                 &1129   & 48.43203 & 51.95         & 8.95e-04  &1293   & 48.41133 & 59.49     & 9.17e-06\\
\hline
MULTADMM            &102    & 48.40617 & 5.19          & 7.03e-04  &113    & 48.40617 & 5.27      & 8.96e-06\\
\hline
PSADMM($m=    2$)   &  744  & 48.42083 & 29.12         & 8.88e-04  & 979   & 48.41131 & 37.57     & 9.16e-06\\
PSADMM($m=    5$)   &  522  & 48.41351 & 15.93         & 9.65e-04  & 883   & 48.41135 & 25.41     & 9.22e-06\\
PSADMM($m=   10$)   &  488  & 48.41123 & 16.46         & 9.49e-04  & 1050  & 48.41128 & 33.57     & 9.49e-06\\
PSADMM($m=   20$)   &  594  & 48.40421 & 25.03         & 9.76e-04  & 1464  & 48.41129 & 58.93     & 9.73e-06\\
PSADMM($m=   50$)   &  986  & 48.40993 & 114.10        & 9.85e-04  & 2000  & 48.41146 & 238.27    & $\sim$\\
PSADMM($m=  100$)   & 1769  & 48.41169 & 250.37        & 9.82e-04  & 2000  & 48.41437 & 270.62    & $\sim$\\
PSADMM($m=  200$)   & 2000  & 48.40564 & 315.62        & $\sim$    & 2000  & 48.40564 & 318.10    & $\sim$\\
\hline
DSADMM($m=    2$)   &  386  & 48.41086 & 39.58         & 9.58e-04  &  623  & 48.41132 & 63.93     & 8.75e-06\\
DSADMM($m=    5$)   &  329  & 48.40631 & 18.11         & 9.68e-04  &  718  & 48.41135 & 40.18     & 9.48e-06\\
DSADMM($m=   10$)   &  412  & 48.41406 & 26.04         & 9.77e-04  &  972  & 48.41130 & 56.93     & 9.14e-06\\
DSADMM($m=   20$)   &  565  & 48.40797 & 42.50         & 9.59e-04  & 1475  & 48.41129 & 111.09    & 9.54e-06\\
DSADMM($m=   50$)   & 1133  & 48.41005 & 149.80        & 9.81e-04  & 2000  & 48.41025 & 261.77    & $\sim$\\
DSADMM($m=  100$)   & 1941  & 48.40482 & 354.21        & $\sim$    & 2000  & 48.41040 & 359.20    & $\sim$\\
DSADMM($m=  200$)   & 2000  & 48.41438 & 369.55        & $\sim$    & 2000  & 48.41438 & 373.04    & $\sim$\\
\hline
\end{tabular}
}
\end{center}

\begin{center}
{\scriptsize\vskip 0.1cm
\begin{tabular}{|l|| c| c| c| c|| c| c| c| c|}
\multicolumn{9}{c}{\small{Table 2.\ Numerical comparison for Basis Pursuit with $n=600$ and $p=2000$}}\\
\hline
   \multicolumn{1}{|c||}{}&\multicolumn{4}{c||}{$\hbox{tol}=$1e-3}&\multicolumn{4}{c|}{$\hbox{tol}=$1e-5} \\
\hline
   \multicolumn{1}{|c||}{\hbox{Method}}&\hbox{Iter}&\hbox{Obj}&\hbox{Time}&\hbox{Error}&\hbox{Iter}&\hbox{Obj}&\hbox{Time}&\hbox{Error} \\
\hline
ADM                 &883    & 96.66843 & 180.72    & 8.56e-04  &1064   & 96.61717 & 217.73    & 9.32e-06\\

\hline
MULTADMM            &66     & 96.65517 & 15.2      & 8.31e-04  &83     & 96.61731 & 20.15     & 9.57e-06\\
\hline
PSADMM($m=    2$)   &  599  & 96.61964 & 112.1     & 9.67e-04  & 845   & 96.61708 & 157.4     & 9.49e-06\\
PSADMM($m=    5$)   &  446  & 96.61319 & 53.6      & 9.75e-04  & 840   & 96.61709 & 100.8     & 9.70e-06\\
PSADMM($m=   10$)   &  445  & 96.61696 & 47.8      & 9.78e-04  & 1060  & 96.61711 & 114.1     & 9.85e-06\\
PSADMM($m=   20$)   &  563  & 96.61880 & 65.7      & 9.87e-04  & 1584  & 96.61710 & 183.6     & 9.77e-06\\
PSADMM($m=   50$)   & 1068  & 96.61684 & 143.9     & 9.84e-04  & 2000  & 96.61755 & 267.5     & $\sim$\\
PSADMM($m=  100$)   & 1916  & 96.61227 & 381.7     & 9.98e-04  & 2000  & 96.60569 & 394.9     & $\sim$\\
PSADMM($m=  200$)   & 2000  & 96.58216 & 612.7     & $\sim$    & 2000  & 96.58216 & 606.9     & $\sim$\\
\hline
DSADMM($m=    2$)   &  757  & 96.65609 & 337.11    & 9.21e-04  & 1000  & 96.61711 & 443.90    & 9.41e-06\\
DSADMM($m=    5$)   &  529  & 96.62273 & 198.85    & 9.78e-04  &  924  & 96.61712 & 346.16    & 9.58e-06\\
DSADMM($m=   10$)   &  502  & 96.61274 & 159.71    & 9.78e-04  & 1109  & 96.61718 & 351.72    & 9.69e-06\\
DSADMM($m=   20$)   &  587  & 96.61563 & 107.28    & 9.75e-04  & 1572  & 96.61720 & 286.33    & 9.74e-06\\
DSADMM($m=   50$)   &  996  & 96.61757 & 228.22    & 9.87e-04  & 2000  & 96.61704 & 457.79    & $\sim$\\
DSADMM($m=  100$)   & 1815  & 96.61811 & 619.64    & 9.96e-04  & 2000  & 96.61395 & 678.21    & $\sim$\\
DSADMM($m=  200$)   & 2000  & 96.59811 & 1001.6    & $\sim$    & 2000  & 96.59811 & 993.54    & $\sim$\\
\hline
\end{tabular}
}
\end{center}

\subsection{Latent Variable Gaussian Graphical Model Selection}
The problem of latent variable Gaussian graphical model selection has been briefly introduced in Section 1. Recall that one of its equivalent reformulation is given by
\begin{eqnarray}
\min_{S,L}&&\langle R, \hat{\Sigma}_X \rangle - \hbox{logdet}(R) + \alpha_1 \| S \|_1 + \alpha_2 {\mathbf{Tr}}(L) + {\mathcal{I}}(L\succeq 0),\label{LVGGMStest}\\
\hbox{s.t.}&&R - S + L = 0\nonumber.
\end{eqnarray}
This model can be viewed as a combination of dimensionality reduction (to identify latent variables) and graphical modeling (to capture remaining statistical structure that is not attributable to
the latent variables). It consistently estimates both the number of hidden components and the conditional graphical model structure among the observed variables. In the following, we show that to solve \eqref{LVGGMStest}, the primal subproblems for {\bf{Algorithm 1}}, {\bf{Algorithm 2}} and {\bf{Algorithm 3}} can be solved exactly and efficiently. To this end, the following two lemmas which can be found in \cite{DaubechiesDefriseDe2004,MaGoldfarbChen2011} are needed.
\begin{lemma}
For $\mu >0$ and $T\in \mathbb{R}^{\ell\times n}$, the solution of the problem
$$
\min_{S\in \mathbb{R}^{\ell\times n}}\ \mu \|S\|_1 + \frac{1}{2}\| S - T \|_F^2,
$$
is given by ${\mathcal{S}}_{\mu}(T)\in \mathbb{R}^{\ell\times n}$, which is defined componentwisely by
\begin{align}
({\mathcal{S}}_{\mu}(T))_{ij}:=\max\{ |T_{ij}| - \mu,0 \}\cdot sgn(T_{ij}).\label{eqMatrixShrink}
\end{align}
\end{lemma}
\begin{lemma}
Let ${\mathcal{S}}_{\mu}$ be defined as in \eqref{eqMatrixShrink}, $K\in \mathbb{R}^{\ell\times n}$ whose rank is $r$ is given, and $\mu > 0$. Then the solution of the problem
$$
\min_{L\in \mathbb{R}^{\ell\times n}}\ \mu \|L\|_* + \frac{1}{2}\| L - K \|_F^2,
$$
is given by ${\mathcal{D}}_{\mu}(K)\in \mathbb{R}^{\ell\times n}$, which is defined by
$$
({\mathcal{D}}_{\mu}(K))_{ij}:=U diag({\mathcal{S}}_{\mu}(\Sigma))V^T,
$$
where $U\in \mathbb{R}^{\ell\times r}$, $V\in \mathbb{R}^{n\times r}$ and $\Sigma \in \mathbb{R}^{r\times r}$ are obtained by the singular value decomposition(SVD) of $K$:
$$
K = U\Sigma V^T\quad \hbox{and}\quad \Sigma = diag(\sigma_1,\sigma_2,\cdots,\sigma_r).
$$
\end{lemma}

Now we are ready to present the steps of different algorithms for solving \eqref{LVGMSeuqal}, by using the previous two lemmas.

 \begin{itemize}
\item {\bf Algorithm 1}: At the $k$th iteration, the update rule is given by:
{\small$$
\qquad R^{k+1} \!=\! U \hbox{diag}(\gamma) U^T, S^{k+1}\!=\!{\mathcal{S}}_{\frac{\alpha_1}{\beta}}\!\!\left(\! R^{k+1} \!+\! L^k \!-\! \frac{\lambda^k}{\beta} \!\right), L^{k+1} \!=\! {\mathcal{D}}_{\frac{\alpha_2}{\beta}}\!\!\left(\! \frac{\lambda^k}{\beta} \!-\! R^{k+1} \!+\! S^{k+1} \!\!\right),
$$}
where $U \hbox{diag}(\sigma)U^T$ is the eigenvalue decomposition of matrix $\frac{1}{\beta}\hat{\Sigma}_{X} \!-\! \frac{1}{\beta}\lambda^k \!-\! S^k \!+\! L^k$ and
$\gamma_i = \left( -\sigma_i + \sqrt{\sigma_i^2 + \frac{4}{\beta}} \right)/2,\ \forall \;i = 1,\cdots,p.$

\item {\bf{Algorithm 2}}: At the $k$th iteration, the update rule is given by:
{\small$$
\qquad R^{k+1}\!=\! U \hbox{diag}(\gamma) U^T,\  S^{k+1}\!=\!{\mathcal{S}}_{\frac{\alpha_1}{\beta}}\!\!\left( - y_2^{k+1} \!-\! \frac{\lambda_2^k}{\beta} \!\right),\ L^{k+1} \!=\! {\mathcal{D}}_{\frac{\alpha_2}{\beta}}\!\!\left( y_3^{k+1} \!+\! \frac{\lambda_3^k}{\beta} \!\right),
$$}
where $U \hbox{diag}(\sigma)U^T$ is the eigenvalue decomposition of another matrix $\frac{1}{\beta}\hat{\Sigma}_{X} \!-\! (y_1^{k+1} \!+\! \frac{1}{\beta}\lambda_1^k)$ and
$\gamma_i = \left( -\sigma_i + \sqrt{\sigma_i^2 + \frac{4}{\beta}} \right)/2,\ \forall\; i = 1,\cdots,p.$

\item {\bf {Algorithm 3}}: At the $k$th iteration, the update rule is given by:
{\small$$
\qquad R^{k+1}\!=\! U \hbox{diag}(\gamma) U^T,\ S^{k+1}\!=\!{\mathcal{S}}_{\frac{\alpha_1}{\beta}}\!\left( - \beta \lambda^{k+1} \!+\! t_2^k \!\right),\ L^{k+1} \!=\! {\mathcal{D}}_{\frac{\alpha_2}{\beta}}\!\left(\! \beta \lambda^{k+1} \!-\! t_3^k \!\right),
$$}
where $U \hbox{diag}(\sigma)U^T$ is the eigenvalue decomposition of another new matrix $\frac{1}{\beta}\hat{\Sigma}_{X} \!-\! (\beta \lambda^{k+1} \!-\! t_1^k)$ and
$\gamma_i = \left( -\sigma_i + \sqrt{\sigma_i^2 + \frac{4}{\beta}} \right)/2,\ \forall\; i = 1,\cdots,p.$
\end{itemize}

In the following, all three methods are compared with PGADM \cite{MaXueZou2012}, which is used to solve the same problem. The stopping criterion is set to be
{\small\begin{eqnarray}
\hbox{Relative Error}\!\!&:=&\!\!\max\!\left\{\! \frac{\|R^{k+1} \!-\! R^k \|_F}{\|R^k\|_F}\!,\!\frac{\|S^{k+1} \!-\! S^k \|_F}{\|S^k\|_F}\!,\!\frac{\|L^{k+1} \!-\! L^k \|_F}{\|L^k\|_F}\!,\right.\nonumber\\
&&\qquad\quad\left. \frac{\|R^{k+1} \!-\! S^{k+1} \!+\! L^{k+1}\|_F}{\max\{ 1,\|R^k\|_F,\|S^k\|_F,\|L^k\|_F \}}\right\} \le \epsilon,\nonumber
\end{eqnarray}}
where $\epsilon$ is some given error tolerance. All the variables are initialized as zero matrices, and the error tolerance $\epsilon$ is set to be $10^{-5}$ in all the  experiments. The comparison results are presented in Table 3, in which `Iter', `Obj' and `Time' denote respectively the iteration number, the objective function value and the CPU time.

\subsubsection{Synthetic Dataset}

We first test the algorithms on synthetic dataset. The sample covariance matrix $\hat{\Sigma}_X$ is generated using the same procedure as in \cite{MaXueZou2012}. Let $p$ and $r$ denote the given dimension of the observed and the latent variables, respectively. We first randomly create a sparse matrix $U \in \mathbb{R}^{(p+r)\times (p+r)}$ with $90\%$ of the entries being zeros, while the nonzero entries were set to be $-1$ or $1$ with equal probability. Then the true covariance matrix is computed by $\Sigma_{X,Y} = (U\cdot U^T)^{-1}$, and the true concentration matrix is given by $\Theta_{X,Y}=U\cdot U^T$. According to $\eqref{eqInverseCovariance}$, the sparse part of ${\Sigma}^{-1}_X$ is given by $$\Theta_X=\Theta_{X,Y}(1:p,1:p),$$
while its low rank part is computed as $\Theta_{XY} \Theta_Y^{-1}\Theta_{YX}=\Theta_{X,Y}(1:p, p+1:p+r)\cdot \Theta_{X,Y}(p+1 : p+r, p+1 : p+r)^{-1}\cdot \Theta_{X,Y}(p+1 : p+r, 1 : p)$.

We then draw $N = 5p$ independent and identically distributed vectors, $Y_1,\cdots ,Y_N$, from the Gaussian distribution ${\mathcal{N}}(0, (\Theta_X - \Theta_{XY} \Theta_Y^{-1}\Theta_{YX})^{-1})$, and compute a sample covariance matrix of  the observed variables according to $\hat{\Sigma}_X:= \frac{1}{N}\sum_{i=1}^N Y_i Y_i^T$. The parameter $\beta$ is set to be $0.1$ for MULTADMM, $0.01$ for PSADMM and DSADMM. For PGADM, following \cite{MaXueZou2012}, $\beta$ is set to be $0.1$, and the parameter $\tau$ in PGADM is set to be $0.6$.

\begin{center}
{\scriptsize\vskip 0.1cm
\begin{tabular}{|c c|| l || c| c| c| c|}
\multicolumn{7}{c}{\small{Table 3. Numerical comparison for LVGGMS}}\\
\hline
   \multicolumn{2}{|c||}{Penalty Parameter}&\multicolumn{1}{c||}{\hbox{Method}}&\hbox{Iter}&\hbox{Obj}&\hbox{Time}&\hbox{Error}\\
\hline
                    &                   &MULTADMM   &50     &-1.831463e+03  &277.45     &9.949149e-06\\
\cline{3-7}
$\alpha_1 = 0.005$  &$\alpha_2 = 0.05$  &PGADM       &162    &-1.835969e+03  &871.57     &9.902557e-06\\
\cline{3-7}
                    &                   &PSADMM &62     &-1.836545e+03  &346.44     &8.954349e-06\\
\cline{3-7}
                    &                   &DSADMM &60     &-1.836545e+03  &340.86     &8.849412e-06\\
\hline
                    &                   &MULTADMM   &19     &-1.739141e+03  &97.41      &7.423706e-06\\
\cline{3-7}
$\alpha_1 = 0.01$       &$\alpha_2= 0.1$        &PGADM       &124    &-1.738876e+03  &626.29     &9.944299e-06\\
\cline{3-7}
                    &                   &PSADMM &47     &-1.739141e+03  &248.06     &9.995064e-06\\
\cline{3-7}
                    &                   &DSADMM &48     &-1.739141e+03  &251.03     &9.957426e-06\\
\hline
                    &                   &MULTADMM   &19     &-1.593774e+03  &67.42      &9.956842e-06\\
\cline{3-7}
$\alpha_1 = 0.02$       &$\alpha_2 = 0.2$       &PGADM       &106    &-1.593676e+03  &534.68     &9.671576e-06\\
\cline{3-7}
                    &                   &PSADMM &37     &-1.593770e+03  &195.42     &8.702799e-06\\
\cline{3-7}
                    &                   &DSADMM &37     &-1.593770e+03  &193.00     &8.719285e-06\\
\hline
                    &                   &MULTADMM   &17     &-1.356307e+03  &89.87      &7.067959e-06\\
\cline{3-7}
$\alpha_1 = 0.04$       &$\alpha_2 = 0.4$       &PGADM       &84     &-1.356285e+03  &420.66     &9.996129e-06\\
\cline{3-7}
                    &                   &PSADMM &35     &-1.356306e+03  &188.52     &9.622188e-06\\
\cline{3-7}
                    &                   &DSADMM &35     &-1.356306e+03  &185.05     &9.622188e-06\\
\hline
\end{tabular}
}
\end{center}

Table 3 compares the algorithms with different choices of penalty parameters $\alpha_1$ and $\alpha_2$. We can find that the proposed methods PSADMM and DSADMM appear to perform better than the state-of-the-art method PGADM: similar objective function values are achieved using significantly less computational time and fewer iterations. Furthermore, MULTADMM is much faster than all the remaining three algorithms, although no theoretical results have been proved yet.

\subsubsection{Stock Dataset}

We then test the algorithms using real stock dataset, which includes the monthly stock return data of companies in the S$\&$P100 index from January 1990 to December 2012. We disregard $26$ companies that were listed after 1990 and only use the data of the remaining $p=74$ companies. The number of months $n$ is equal to 276. Each correlation coefficient between two stocks is computed based on the $n$ dimensional S$\&$P100 index vector of each stock. As a result the estimated covariance matrix is a $p\times p$ matrix with diagonal elements all being $1$. We choose the penalty parameters as $\alpha_1=0.005$ and $\alpha_2=0.01$; choose $\beta = 10$ for all four methods; choose $\tau = 0.6$ for PGADM. The performance of the algorithm is presented in Table 4. The identified relationship among the companies, characterized by the sparse part of the inverse of the estimated concentration matrix $\hat{\Sigma}_X$, is shown in Figure 2.

\begin{center}
{\scriptsize\vskip 0.1cm
\begin{tabular}{| l || c| c| c| c|}
\multicolumn{5}{c}{\small{Table 4. Numerical comparison for S$\&$P100 index}}\\
\hline
   \multicolumn{1}{|c||}{\hbox{Method}}&\hbox{Iter}&\hbox{Obj}&\hbox{Time}&\hbox{Error}\\
\hline
MULTADMM    &193    &-48.55483  &1.01   &9.502830e-06\\
\hline
PGADM        &393    &-48.45214  &1.97   &9.209239e-06\\
\hline
PSADMM      &380    &-48.57557  &1.90   &9.555332e-06\\
\hline
DSADMM      &373    &-48.57557  &1.86   &9.356121e-06\\
\hline
\end{tabular}
}
\end{center}

\begin{figure}[h]
\centering
\includegraphics[scale=0.55]{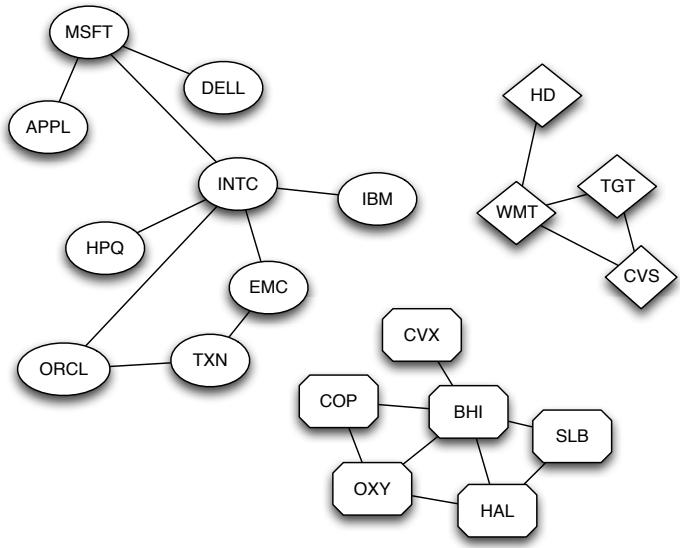}
\caption{Relationship performance between some companies}
\end{figure}

From Table 4, we see that MULTADMM still performs the best with substantially less computational time and fewer number of iterations, while the PSADMM and DSADMM slightly outperform PGADM. Figure 2 shows the identified relationship between some of those $74$ chosen companies in S$\&$P100. Recall that two companies $i,j$ are related if the $(i,j)$th entry of the computed sparse part $S$ is nonzero. All these methods are applied to the same data-set, as a result all the results are almost the same. So that we can choose any one of the obtained $S$ to plot Figure 2, and the MULTADMM result is chosen in this paper. In Figure 2, three groups of the companies are shown, and it is clear that the identified groups are intuitively meaningful: they represent {\em{Information Technology Companies}}, {\em{Retail Companies}} and {\em{Oil and Gas Companies}} respectively.

\subsection{Robust Principal Component Analysis}

In this experiment, we test the proposed methods on the RPCA model (\ref{RobustPCA}), applied to the video background extraction problem with missing and noisy data. It is well known that a segment of video, which consists of several consecutive frames, is a good candidate for sparse and low-rank modeling. In particular, let $M$ denote the matrix representation of a sequence of video frames where each column represents a frame. Then the matrix representation of the background should have low rank for the reason that the background needs to be flexible enough to accommodate changes in the scene. On the other hand, the foreground objects such as cars or pedestrians can be collected by a sparse matrix with only a small number of nonzero elements.
Moreover, in real application the video may include noisy and even missing pixels. Suppose we have a sequence of video
consisting of $n$ frames with each frame having $\ell$ pixels. Mathematically, let us use $L\in \mathbb{R}^{\ell\times n}$, $S\in \mathbb{R}^{\ell\times n}$ and $Z\in \mathbb{R}^{\ell\times n}$ to denote the background, foreground and the noisy pixels. Let the index set $\Omega$ indicate the locations of observed pixels, i.e., pixels outside $\Omega$ are missing. Then the background extraction problem for the noise corrupted video can be stated as follows
\begin{eqnarray}
\min_{L,S,Z} && \|L\|_* + \tau \|S\|_1,\label{problemRPCAVideo}\\
\hbox{s.t.} && L + S + Z = M,\nonumber\\
&&Z\in {\mathcal{H}}:=\{Z\in \mathbb{R}^{\ell\times n}|\|P_{\Omega}(Z)\|_F\le \delta\}.\nonumber
\end{eqnarray}
For detailed discussion about this model, we refer the interested readers to \cite{CandesRPCA2011}.

Next we specialize {\bf{Algorithm 1}}, {\bf{Algorithm 2}} and {\bf{Algorithm 3}} to solve problem \eqref{problemRPCAVideo}.

For {\bf{Algorithm 1}}, at the $k$th iteration, the variables are updated according to
\begin{eqnarray}\label{subAlg1RPCA}
L^{k+1} \!\!\!\!&=&\!\!\!\! {\mathcal{D}}_{\frac{1}{\beta}}\left( M + \frac{\lambda^k}{\beta} - S^k - Z^k \right),\qquad S^{k+1} = {\mathcal{S}}_{\frac{\tau}{\beta}}\left( M + \frac{\lambda^k}{\beta} - L^{k+1} - Z^k \right),\nonumber\\
&& (Z^{k+1})_{ij}= \left\{\begin{array}{l}
                        N^k_{ij},\qquad\qquad\qquad\quad\  \hbox{if} (i,j)\ \notin\Omega,\\
                        \frac{\min\{\|P_{\Omega}(N^k)\|_F,\delta\}}{\|P_{\Omega}(N^k)\|_F} N_{ij}^k,\ \hbox{if}\  (i,j)\in \Omega, \\
                        \end{array}\right.\nonumber
\end{eqnarray}
where $N^k: = M + \frac{\lambda^k}{\beta} - L^{k+1} - S^{k+1}$.

For {\bf{Algorithm 2}}, at $k$th iteration, the variables are updated according to
\begin{eqnarray}\label{subAlg2RPCA}
L^{k+1} \!\!\!\!&=&\!\!\!\! {\mathcal{D}}_{\frac{1}{\beta}}\left( \frac{M}{3} + y_1^{k+1} + \frac{\lambda_1^k}{\beta} \right),\qquad S^{k+1} = {\mathcal{S}}_{\frac{\tau}{\beta}}\left( \frac{M}{3} + y_2^{k+1} + \frac{\lambda_2^k}{\beta} \right),\nonumber\\
&&(Z^{k+1})_{ij}= \left\{\begin{array}{l}
                        N^k_{ij},\qquad\qquad\qquad\quad\  \hbox{if} (i,j)\ \notin\Omega,\\
                        \frac{\min\{\|P_{\Omega}(N^k)\|_F,\delta\}}{\|P_{\Omega}(N^k)\|_F} N_{ij}^k,\ \hbox{if}\  (i,j)\in \Omega, \\
                        \end{array}\right.\nonumber
\end{eqnarray}
where $y_i$ is the slack variables defined {\bf{Algorithm 2}}; $N^k := \frac{M}{3} + y_3^{k+1} + \frac{\lambda_3^k}{\beta}$.

As for {\bf{Algorithm 3}}, the variables are updated according to
\begin{eqnarray}\label{subAlg3RPCA}
L^{k+1} \!\!\!\!&=&\!\!\!\! {\mathcal{D}}_{\beta}\left( \beta \lambda^{k+1} - t_1^k \right),\qquad S^{k+1}\ =\ {\mathcal{S}}_{\tau\beta}\left( \beta \lambda^{k+1} - t_2^k \right),\nonumber\\
&&(Z^{k+1})_{ij}= \left\{\begin{array}{l}
                        N^k_{ij},\qquad\qquad\qquad\quad\  \hbox{if} (i,j)\ \notin\Omega,\\
                        \frac{\min\{\|P_{\Omega}(N^k)\|_F,\delta\}}{\|P_{\Omega}(N^k)\|_F} N_{ij}^k,\ \hbox{if}\  (i,j)\in \Omega, \\
                        \end{array}\right.\nonumber
\end{eqnarray}
where $t_i$ is the slack variable defined in {\bf{Algorithm 3}}; $N^k := \beta \lambda^{k+1} - t_3^k$.

The test dataset is a video taken at the hall of an aiport \footnote{Available at http://perception.i2r.a-star.edu.sg/$\hbox{bk}\_\hbox{model}$/$\hbox{bk}\_\hbox{index}$.html.} that consists of 200 grayscale frames, each of the size 144$\times$176. As a result the matrix $M$ is of dimension 25,344$\times$200. The index set $\Omega$ is determined randomly with a fixed sampling ratio $sr=80\%$, meaning that $20\%$ of the pixels are missing. The proposed three methods (MULTADMM, PSADMM and DSADMM) are compared with ADMMGBS in \cite{HeTaoYuan2012} and VASALM in \cite{TaoYuan2011}. The parameters in model (\ref{RobustPCA}) is set as $\tau = \frac{1}{\sqrt{l}}$ and $\delta = 10^{-2}$. We choose $\beta = \frac{0.002}{\|M\|_1}$ for MULTADMM and VASALM, $\beta = \frac{0.004}{\|M\|_1}$ for PSADMM, and $\beta = 2\cdot\|M\|_1$ for DSADMM. All the algorithms are initialized using zero matrices. The stopping criterions are set to be
$$
\hbox{Relative Error}:=\max \left\{\frac{\| L^{k+1} - L^k \|_F}{1 + \|L^k\|_F},\frac{\| S^{k+1} - S^k \|_F}{1 + \|S^k\|_F}, \frac{|f^{k+1} - f^k|}{|f^k|}\right\}\le 10^{-3},
$$
where $\left\{f^k := \|L^k\|_* + \tau \|S^k\|_1|\right\}$ is the sequence of objective value of \eqref{MultipleBlockProblem}.
In Table 5, we have defined $\hbox{Obj}_{error}$ as the relative successive difference of the objective, i.e.,
$$\hbox{Obj}_{error} := \frac{|f^{k+1} - f^k|}{|f^k|}.$$
Moreover, ``Iter", ``Time", ``nnz(S)", and ``rank(L)" respectively denote the iteration number, the computation time, the number of nonzero elements in $S$, and the rank of $L$.

\begin{center}
{\scriptsize\vskip 0.1cm
\begin{tabular}{| l || c| c| c| c| c|}
\multicolumn{6}{c}{\small{Table 5. Numerical comparison for RPCA}}\\
\hline
\hbox{Method}&\hbox{Iter}&\hbox{Time}&\hbox{rank(L)}&\hbox{nnz(S)}&$\hbox{Obj}_{error}$\\
\hline
VASALM      &24 &42.84  &  4    & 149178    &5.72e-04\\
\hline
ADMGBS      &25 &44.70  &  4    & 149656    &9.14e-05\\
\hline
MULTADMM    &22 &40.66  &  4    & 146708    &3.27e-04\\
\hline
PSADMM      &25 &43.88  &  4    & 148015    &6.41e-04\\
\hline
DSADMM      &26 &45.27  &  4    & 150842    &3.72e-03\\
\hline
\end{tabular}
}
\end{center}

From the results in Table 5, once again we observe that the MULTADMM is the most efficient algorithm for solving RPCA, although the margin of advantage over the remaining algorithms is small.

\section{Conclusion}\label{secConclusion}

In this paper, we systematically study several ADMM based algorithms for solving multiple block separable convex minimization problems. Two algorithms are proposed, which apply the classical two-block ADMM to either a primal or a dual reformulation of the original multi-block problem. Various theoretical properties of these two algorithms, such as their global convergence, iteration complexity and linear convergence, are shown by leveraging existing results. Through extensive numerical experiments, we observe that these two methods are competitive with state-of-the-art algorithms. Somewhat surprisingly, we show that in most cases, the classical multiple block ADMM method is computationally much more efficient than the proposed algorithms. This observation suggests that it is important to investigate the theoretical properties of the multiple block ADMM methods.

\medskip
{\bf{Acknowledgment}}: The authors are grateful to Professor Bingsheng He of Nanjing University and Dr.\ Tsung-Hui Chang of National Taiwan University of Science and Technology for valuable suggestions and discussions.

\section{Appendix}\label{seeAppendix}
{\bf{Proof of Theorem \ref{ComplexityObjective}}}. First recall that at ($k+1$)-th iteration, the optimality conditions of each subproblem about $x_i$ for both exact and inexact versions are given by
\begin{align}
\left(x_i - x_i^{k+1}\right)^T \!\!\left\{ \partial f_i (x_i^{k+1}) + \beta A_i^T \left( A_i x_i^{k+1} \!-\! \frac{b}{m} - y_i^{k+1} - \frac{\lambda_i^k}{\beta} \right) \!+\! P_i \left( x_i^{k+1} - x_i^k \right) \right\}\!\!\ge 0&,\nonumber\\
\vspace{-4cm}\forall\ x_i\in {\mathcal{X}}_i&\nonumber
\end{align}
where $P_i = 0$ for exact version and $P_i = \tau_i I - \beta A_i^T A_i$ for inexact version. Further substitute $\lambda_i^{k+1} = \lambda_i^k - \beta\left( A_i x_i^{k+1} - \frac{b}{m} - y_i^{k+1} \right)$ into the above inequalities,
$$
\left(x_i - x_i^{k+1}\right)^T \left\{ \partial f_i (x_i^{k+1}) - A_i^T \lambda_i^{k+1} + P_i \left( x_i^{k+1} - x_i^k \right) \right\}\ge 0,\ \forall\ x_i\in {\mathcal{X}}_i,
$$
then combine with the convexity of $f_i(\cdot)$ and let $x_i = x_i^*\in {\mathcal{X}}_i$
{\small\begin{eqnarray}\label{AppendixBasic}
\!\!\!&&\!\!\!\!\!f_i(x_i^{k+1}) - f_i(x_i^*)\le \left( A_i x_i^{k+1} - A_i x_i^* \right)^T \lambda_i^{k+1} + \left( x_i^* - x_i^{k+1} \right)^T P_i \left( x_i^{k+1} - x_i^k \right)\nonumber\\
\!\!\!&=&\!\!\!\!\!\! \left( y_i^{k+1} - y_i^* + \frac{1}{\beta} \left( \lambda_i^k - \lambda_i^{k+1} \right) \right)^T \lambda_i^{k+1} + \left( x_i^* - x_i^{k+1} \right)^T P_i \left( x_i^{k+1} - x_i^k \right)\nonumber\\
\!\!\!&=&\!\!\!\!\!\! \left( y_i^{k+1} \!-\! y_i^* \right)^T\!\! \lambda_i^{k+1} \!\!+\! \frac{\left\|\lambda_i^k\right\|^2 \!\!\!-\! \left\| \lambda_i^{k+1} \right\|^2 \!\!\!-\! \left\| \lambda_i^k \!-\! \lambda_i^{k+1} \right\|^2}{2\beta} \!+\! \frac{\left\|x_i^* \!\!-\! x_i^k \right\|_{P_i}^2 \!\!\!-\! \left\| x_i^* \!\!-\! x_i^{k+1} \right\|_{P_i}^2}{2}.
\end{eqnarray}}
At $(k+1)$-th iteration,  the optimality condition of the subproblem about $\bfy$ is given by
$$
\sum_{i=1}^m \left( y_i - y_i^{k+1} \right)^T \left( y_i^{k+1} - \left( A_i x_i^k - \frac{b}{m} - \frac{\lambda_i^k}{\beta}\right) \right) \ge 0,\ \forall\ \bfy\in {\mathcal{Y}},
$$
\begin{equation}\label{temp1}
\Rightarrow\qquad \sum_{i=1}^m \left( y_i - y_i^{k+1} \right)^T \left( \lambda_i^{k+1} + \beta A_i x_i^{k+1} - \beta A_i x_i^k \right) \ge 0,\ \forall\ \bfy\in {\mathcal{Y}},
\end{equation}
after letting $\bfy = \bfy^*\in {\mathcal{Y}}$ in \eqref{temp1}, we obtain
{\footnotesize\begin{eqnarray}
&&\sum_{i=1}^m \left( y_i^{k+1} - y_i^* \right)^T \lambda_i^{k+1} \le \beta \sum_{i=1}^m \left( y_i^* - y_i^{k+1} \right)^T \left[ \left( A_i x_i^{k+1} - \frac{b}{m}\right)  - \left( A_i x_i^k - \frac{b}{m}\right) \right]\nonumber\\
\!\!\!\!&=&\!\!\!\! \frac{\beta}{2} \sum_{i=1}^m \left\| y_i^* \!-\! A_i x_i^k \!+\! \frac{b}{m} \right\|^2 \!-\! \left\| y_i^* \!-\! A_i x_i^{k+1} \!+\! \frac{b}{m} \right\|^2 \!+\! \left\| y_i^{k+1} \!-\! A_i x_i^{k+1} \!+\! \frac{b}{m} \right\|^2 \!-\! \left\| y_i^{k+1} \!-\! A_i x_i^k \!+\! \frac{b}{m} \right\|^2\nonumber\\
\!\!\!\!&\le&\!\!\!\! \frac{\beta}{2} \sum_{i=1}^m \left\{ \left\| y_i^* \!-\! A_i x_i^k \!+\! \frac{b}{m} \right\|^2 \!-\! \left\| y_i^* \!-\! A_i x_i^{k+1} \!+\! \frac{b}{m} \right\|^2 \!+\! \left\| y_i^{k+1} \!-\! A_i x_i^{k+1} \!+\! \frac{b}{m} \right\|^2\right\}.\nonumber
\end{eqnarray}}
Summing \eqref{AppendixBasic} from $i = 1$ to $i = m$, and using the above inequality, we can bound the difference of the current and optimal objective value by
{\footnotesize\begin{eqnarray}
&&\theta(\bfx^{k+1}) - \theta(\bfx^*)\le \sum_{i=1}^m \frac{\left\|\lambda_i^k\right\|^2 \!-\! \left\| \lambda_i^{k+1} \right\|^2 \!-\! \left\| \lambda_i^k \!-\! \lambda_i^{k+1} \right\|^2}{2\beta} \!+\! \frac{\left\|x_i^* \!-\! x_i^k \right\|_{P_i}^2 \!-\! \left\| x_i^* \!-\! x_i^{k+1} \right\|_{P_i}^2}{2}\nonumber\\
&&\qquad +\  \frac{\beta}{2} \sum_{i=1}^m \left\{ \left\| y_i^* \!-\! A_i x_i^k \!+\! \frac{b}{m} \right\|^2 \!-\! \left\| y_i^* \!-\! A_i x_i^{k+1} \!+\! \frac{b}{m} \right\|^2 \!+\! \left\| y_i^{k+1} \!-\! A_i x_i^{k+1} \!+\! \frac{b}{m} \right\|^2\right\}\nonumber\\
&\le&\!\!\!\!\!\!\! \sum_{i=1}^m\! \left\{\!\! \frac{\left\|\lambda_i^k\right\|^2 \!\!\!-\! \left\| \lambda_i^{k+1} \right\|^2}{2\beta} \!+\! \frac{\left\|x_i^* \!\!-\! x_i^k \right\|_{P_i}^2 \!\!\!-\! \left\| x_i^* \!-\! x_i^{k+1} \right\|_{P_i}^2}{2} \!+\! \beta\frac{\left\| y_i^* \!-\! A_i x_i^k \!+\! \frac{b}{m} \right\|^2 \!-\! \left\| y_i^* \!-\! A_i x_i^{k+1} \!+\! \frac{b}{m} \right\|^2}{2} \right\},\nonumber
\end{eqnarray}}
while the last inequality is obtained from the multiplier update step
$$y_i^{k+1} - A_i x_i^{k+1} + \frac{b}{m} = \frac{1}{\beta}(\lambda_i^{k+1} - \lambda_i^k).$$
Finally define $\tilde{\bfx}^K = \frac{1}{K} \sum_{k=1}^K \bfx^k$,
{\footnotesize\begin{eqnarray}
&&\!\!\!\!\theta(\tilde{\bfx}^K) - \theta(\bfx^*) \le \frac{1}{K} \sum_{k=1}^K \left( \theta(\bfx^k) - \theta(\bfx^*) \right)\nonumber\\
&\le&\!\!\!\! \frac{1}{K} \sum_{i=1}^m \left\{ \frac{\left\|\lambda_i^0\right\|^2 \!\!-\! \left\| \lambda_i^K \right\|^2}{2\beta} \!+\! \frac{\left\|x_i^* \!-\! x_i^0 \right\|_{P_i}^2 \!-\! \left\| x_i^* \!-\! x_i^K \right\|_{P_i}^2}{2} \!+\! \beta\frac{\left\| y_i^* \!-\! A_i x_i^0 \!+\! \frac{b}{m} \right\|^2 \!-\! \left\| y_i^* \!-\! A_i x_i^K \!+\! \frac{b}{m} \right\|^2}{2} \right\}\nonumber\\
&\le&\!\!\!\! \frac{1}{K} \sum_{i=1}^m \left\{ \frac{\|\lambda_i^0\|^2}{2\beta} + \frac{\|x_i^* - x_i^0\|^2}{2} + \frac{\beta \| y_i^* - A_i x_i^0 + \frac{b}{m} \|^2}{2} \right\} \thicksim \mathcal{O}\left(\frac{1}{K}\right).\qquad\qquad\qquad\qquad\qquad\qquad\qquad\square\nonumber
\end{eqnarray}}

{\bf{Proof of Theorem \ref{Algorithm2Convergence} for inexact version}}. First recall the optimality conditions of each subproblems about $x_i$ and $\lambda$ in the ($k+1$)-th iteration, which are expressed as $\forall\ x_i\in {\mathcal{X}}_i$
$$
\left(x_i - x_i^{k+1}\right)^T \!\!\left\{ \partial f_i (x_i^{k+1}) + \tau_i \left[ x_i^{k+1} - x_i^k + \frac{1}{\tau_i \beta} A_i^T \left( A_i x_i^k - \beta \lambda^{k+1} - t_i^k \right) \right] \right\}\ge 0,
$$
$$
\left( \lambda - \lambda^{k+1} \right)^T \left\{ -b - \sum\limits_{i=1}^m t_i^k + \beta \sum\limits_{i=1}^m \left( \lambda^{k+1} - \lambda_i^k \right) \right\}\ge 0,\ \forall\  \lambda,
$$
where $\lambda_i$ and $t_i$ also satisfy the following two equalities,
\begin{align}\label{eq:AppEquality}
\lambda_i^{k+1} = \lambda^{k+1} - \frac{1}{\beta} \left( t_i^k + A_i x_i^{k+1}\right),\qquad \hbox{and}\qquad t_i^{k+1} = t_i^k - \beta \left(\lambda^{k+1} - \lambda_i^{k+1}\right).
\end{align}
After substituting these two equalities into the above optimality conditions, we can obtain
\begin{equation}\label{DualSplitOptimalityx}
\left(x_i - x_i^{k+1}\right)^T \left\{ \partial f_i (x_i^{k+1}) - A_i^T \lambda^{k+1} + \tau_i \left( x_i^{k+1} - x_i^k \right) \right\}\ge 0,\ \forall\ x_i\in {\mathcal{X}}_i,
\end{equation}
\begin{equation}\label{DualSplitOptimalitylambda}
\left( \lambda - \lambda^{k+1} \right)^T \left\{ \sum\limits_{i=1}^m A_i x_i^{k+1} - b + \beta \sum\limits_{i=1}^m \left( \lambda_i^{k+1} - \lambda_i^k \right) \right\}\ge 0,\ \forall\  \lambda.
\end{equation}
Adding these inequalities together, we have
\begin{eqnarray}\label{AppendixOptimality}
\qquad&&\!\!\!\!\sum_{i=1}^m \left(x_i - x_i^{k+1}\right)^T \left( \partial f_i (x_i^{k+1}) - A_i^T \lambda^{k+1} \right) + \left( \lambda - \lambda^{k+1} \right)^T \left( \sum\limits_{i=1}^m A_i x_i^{k+1} \!-\! b \right)\nonumber\\
\qquad&\ge&\!\!\!\! \sum_{i=1}^m \left\{ \tau_i \left(x_i - x_i^{k+1}\right)^T \left( x_i^k - x_i^{k+1} \right) + \beta \left( \lambda - \lambda^{k+1} \right)^T \left( \lambda_i^k - \lambda_i^{k+1} \right) \right\}.
\end{eqnarray}
Set $x_i = x_i^*$ and $\lambda = \lambda^*$ where ($x_1^*,\cdots,x_m^*,\lambda^*$) is any saddle point of the Lagrangian function of \eqref{MultipleBlockProblem}, then clearly the left hand side of the above inequality becomes no greater than zero. Note that the following equality is true
\begin{eqnarray}
\left( \lambda^* \!-\! \lambda^{k+1} \right)^T\!\! \left( \lambda_i^k \!-\! \lambda_i^{k+1} \right) \!\!\!\!\!&=&\!\!\!\!\! \left( \lambda_i^* - \lambda_i^{k+1} \right)^T\!\! \left( \lambda_i^k - \lambda_i^{k+1} \right) + \left( \lambda_i^{k+1} - \lambda^{k+1} \right)^T\!\! \left( \lambda_i^k - \lambda_i^{k+1} \right)\nonumber\\
\!\!\!\!\!&=&\!\!\!\!\! \left( \lambda_i^* \!-\! \lambda_i^{k+1} \right)^T\!\! \left( \lambda_i^k \!-\! \lambda_i^{k+1} \right) \!+\! \left[ A_i \left( x_i^k \!-\! x_i^{k+1} \right) \right]^T\!\! \left( \lambda_i^k \!-\! \lambda_i^{k+1} \right).\nonumber
\end{eqnarray}
Then we can rewrite \eqref{AppendixOptimality} by
\begin{eqnarray}
&&\!\!\!\!\sum_{i=1}^m \left\{ \tau_i \left(x_i^* - x_i^k\right)^T \left( x_i^k - x_i^{k+1} \right) + \beta \left( \lambda_i^* - \lambda_i^k \right)^T \left( \lambda_i^k - \lambda_i^{k+1} \right) \right\}\nonumber\\
&\le&\!\!\!\! - \sum_{i=1}^m \left\{ \tau_i \left\|x_i^k - x_i^{k+1}\right\|^2 + \beta \left\| \lambda_i^k - \lambda_i^{k+1} \right\|^2 + \left[ A_i \left( x_i^k - x_i^{k+1} \right) \right]^T \left( \lambda_i^k - \lambda_i^{k+1} \right) \right\}, \nonumber
\end{eqnarray}
Let us define a matrix ${\mathcal{T}} = \hbox{diag}\{ \tau_1 I,\cdots,\tau_m I \}\in {\mathbb{R}}^{n\times n}$. Using the previous inequality,  we have
{\footnotesize\begin{eqnarray}
&&\left\|\bfx^{k+1} \!- \bfx^*\right\|_{\mathcal{T}}^2 \!+\! \beta \left\| \tilde{\bflambda}^{k+1} \!- \tilde{\bflambda}^* \right\|^2 = \left\|\bfx^{k+1} \!- \bfx^k \!+ \bfx^k \!- \bfx^*\right\|_{\mathcal{T}}^2\! +\! \beta \left\| \tilde{\bflambda}^{k+1} \!- \tilde{\bflambda}^k \!+ \tilde{\bflambda}^k \!- \tilde{\bflambda}^* \right\|^2\nonumber\\
\!\!&=&\!\! \left\|\bfx^k \!- \bfx^*\right\|_{\mathcal{T}}^2 \!+\! \beta \left\| \tilde{\bflambda}^k \!- \tilde{\bflambda}^* \right\|^2 \!+\! 2 \left( \bfx^k \!- \bfx^* \right)^T {\mathcal{T}} \left( \bfx^{k+1} \!- \bfx^k \right)\! + \!2 \beta \left( \tilde{\bflambda}^k \!- \tilde{\bflambda}^* \right)^T\!\! \left( \tilde{\bflambda}^{k+1} \!- \tilde{\bflambda}^k \right)\nonumber\\
&&\!+\!\left\|\bfx^k \!- \bfx^{k+1}\right\|_{\mathcal{T}}^2 \!+\! \beta \left\| \tilde{\bflambda}^k \!- \tilde{\bflambda}^{k+1} \right\|^2 \nonumber\\
\!\!&\le&\!\! \left\|\bfx^k \!- \bfx^*\right\|_{\mathcal{T}}^2 \!+\! \beta \left\| \tilde{\bflambda}^k \!- \tilde{\bflambda}^* \right\|^2 \!-\!\left\|\bfx^k \!- \bfx^{k+1}\right\|_{\mathcal{T}}^2 \!-\! \beta \left\| \tilde{\bflambda}^k \!- \tilde{\bflambda}^{k+1} \right\|^2 - 2\sum_{i=1}^m \left[ A_i \left( x_i^k \!- x_i^{k+1} \right) \right]^T \!\!\nonumber\\
&&  \left( \lambda_i^k \!- \lambda_i^{k+1} \right), \nonumber\\
\!\!&\le&\!\! \left\|\bfx^k \!- \bfx^*\right\|_{\mathcal{T}}^2 \!+\! \beta \left\| \tilde{\bflambda}^k \!- \tilde{\bflambda}^* \right\|^2 \!-\sum_{i=1}^m \left\| \left( \begin{array}{c} x_i^k\\ \lambda_i^k \end{array} \right) - \left( \begin{array}{c} x_i^{k+1}\\ \lambda_i^{k+1} \end{array} \right) \right\|_{{\mathcal{Q}}_i} \nonumber
\end{eqnarray}}
when $ \tilde{\bflambda} = \left( \lambda_1^T,\cdots,\lambda_m^T \right)^T $ and ${\mathcal{Q}}_i = \left( \begin{array}{cc}\tau_i I\!&\!A_i^T\\ A_i\!&\!\beta I\end{array} \right)\succ 0$ if $\tau_i > \frac{1}{\beta} \rho (A_i^T A_i) $. In summary $\{\left\|\bfx^k \!- \bfx^*\right\|_{\mathcal{T}}^2 \!+\! \beta \left\| \tilde{\bflambda}^k \!- \tilde{\bflambda}^* \right\|^2\}$ is a non-increasing sequence satisfying
{\small$$
\left\|\bfx^{k+1} \!- \bfx^*\right\|_{\mathcal{T}}^2 \!+\! \beta \left\| \tilde{\bflambda}^{k+1} \!- \tilde{\bflambda}^* \right\|^2 \le \left\|\bfx^k \!- \bfx^*\right\|_{\mathcal{T}}^2 \!+\! \beta \left\| \tilde{\bflambda}^k \!- \tilde{\bflambda}^* \right\|^2 - \left\| \left( \begin{array}{c} x^k\\ \tilde{\lambda}^k \end{array} \right) - \left( \begin{array}{c} x^{k+1}\\ \tilde{\lambda}^{k+1} \end{array} \right) \right\|_{{\mathcal{Q}}}.
$$}
where ${\mathcal{Q}} = \hbox{diag}({\mathcal{Q}}_1,\cdots,{\mathcal{Q}}_m)$. We can get the following result by summing from $1$ to $\infty$,
$$
\sum_{i=1}^{\infty} \left\| \left( \begin{array}{c} x^k\\ \tilde{\lambda}^k \end{array} \right) - \left( \begin{array}{c} x^{k+1}\\ \tilde{\lambda}^{k+1} \end{array} \right) \right\|_{{\mathcal{Q}}} \le \left\|\bfx^0 - \bfx^*\right\|_{\mathcal{T}}^2 + \beta \left\| \tilde{\bflambda}^0 - \tilde{\bflambda}^* \right\|^2 < \infty,
$$
which implies
$$
\left\| \bfx^k - \bfx^{k+1} \right\|^2 \rightarrow 0,\qquad\hbox{and}\qquad \left\| \tilde{\bflambda}^k - \tilde{\bflambda}^{k+1} \right\|^2 \rightarrow 0.
$$
Using the equality relationship \eqref{eq:AppEquality} we can also obtain,
$$
\left\| \bft^k - \bft^{k+1} \right\|^2 \rightarrow 0,\qquad\hbox{and}\qquad \left\| \lambda^k - \lambda^{k+1} \right\|^2 \rightarrow 0.
$$
Let $\bfx^{k_j}$ and $\lambda^{k_j}$ be the subsequence of $\bfx^k$ and $\lambda^k$ which converge to $\bfx^{\infty}$ and $\lambda^{\infty}$. Then recall the optimality conditions of this algorithm, $\bfx^{\infty}$ and $\lambda^{\infty}$ can be easily proved to be the primal and dual optimal solution of \eqref{MultipleBlockProblem} respectively. Further with the non-increasing property of $\{\|\bfx^k \!- \bfx^{\infty}\|_{\mathcal{T}}^2 \!+\! \beta \| \tilde{\bflambda}^k \!- \tilde{\bflambda}^{\infty} \|^2\}$, $\{\bfx^k\}$ and $\lambda^k$ must converge to $\bfx^{\infty}$ and $\lambda^{\infty}$. From \eqref{eq:AppEquality} we have $t_i^k = - A_i x_i^k$, which implies the convergence of $\{\tilde{\bflambda}^k\}$ and $\{\bft^k\}$. \hfill$\square$
\bigskip\\
{\bf{Proof of Corollary \ref{corComplexity}}}. Let us rewrite the optimality conditions of each subproblems in {\bf{Algorithm~3}},
$$
\left(x_i - x_i^{k+1}\right)^T \left\{ \partial f_i (x_i^{k+1}) - A_i^T \lambda^{k+1} + G_i \left( x_i^{k+1} - x_i^k \right) \right\}\ge 0,\ \forall\ x_i\in {\mathcal{X}}_i,
$$
$$
\left( \lambda - \lambda^{k+1} \right)^T \left\{ \sum\limits_{i=1}^m A_i x_i^{k+1} - b + \beta \sum\limits_{i=1}^m \left( \lambda_i^{k+1} - \lambda_i^k \right) \right\}\ge 0,\ \forall\  \lambda.
$$
where $G_i = \frac{1}{\beta} A_i^T A_i$ for exact version and $G_i = \tau_i I$ for inexact version. Using the Lagrangian function of \eqref{MultipleBlockProblem} and for all $(x_1,\cdots,x_m,\lambda)$, we obtain
{\small\begin{eqnarray}
&& L(x_1^{k+1},\cdots,x_m^{k+1},\lambda) - L(x_1,\cdots,x_m,\lambda^{k+1})\nonumber\\
&=& \sum\limits_{i=1}^m f_i (x_i^{k+1}) - \sum\limits_{i=1}^m f_i (x_i) -
\langle \lambda, \sum\limits_{i=1}^m A_i x_i^{k+1} - b \rangle + \langle \lambda^{k+1}, \sum\limits_{i=1}^m A_i x_i - b \rangle\nonumber\\
&\le & \sum\limits_{i=1}^m \left(x_i^{k+1} - x_i\right)^T \left\{ \partial f_i (x_i^{k+1}) - A_i^T \lambda^{k+1} \right\} + \left( \lambda^{k+1} - \lambda \right)^T \left\{ \sum\limits_{i=1}^m A_i x_i^{k+1} - b \right\}\nonumber\\
&\le & \sum\limits_{i=1}^m \left(x_i - x_i^{k+1}\right)^T G_i \left( x_i^{k+1} - x_i^k \right) + \beta \sum\limits_{i=1}^m \left( \lambda - \lambda^{k+1} \right)^T \left( \lambda_i^{k+1} - \lambda_i^k \right)\nonumber\\
& = & \sum\limits_{i=1}^m \left\{ \frac{1}{2} \left( \| x_i \!-\! x_i^k \|_{G_i}^2 \!-\! \| x_i \!-\! x_i^{k+1} \|_{G_i}^2 \!-\! \| x_i^{k+1} \!-\! x_i^k \|_{G_i}^2 \right) \!+\! \frac{\beta}{2} \left( \| \lambda \!-\! \lambda_i^k \|^2 \!-\! \| \lambda \!-\! \lambda_i^{k+1} \|^2 \right.\right.\nonumber\\
&& \left.\left.+ \| \lambda^{k+1} - \lambda_i^{k+1} \|^2 - \| \lambda^{k+1} - \lambda_i^k \|^2 \right) \right\}.\nonumber
\end{eqnarray}}
Next we bound part of the right hand side in the above inequality. Notice that
\begin{eqnarray}
&& \| x_i^{k+1} - x_i^k \|_{G_i}^2 + \beta \left( \| \lambda^{k+1} - \lambda_i^k \|^2 - \| \lambda^{k+1} - \lambda_i^{k+1} \|^2 \right)\nonumber\\
&=& \| x_i^{k+1} - x_i^k \|_{G_i}^2 + 2\beta \left( \lambda^{k+1} - \lambda_i^{k+1} \right)^T \left( \lambda_i^{k+1} - \lambda_i^k \right) + \beta \| \lambda_i^{k+1} - \lambda_i^k \|^2\nonumber\\
&=& \| x_i^{k+1} - x_i^k \|_{G_i}^2 + 2\beta \left( x_i^{k+1} - x_i^k \right)^T A_i^T \left( \lambda_i^{k+1} - \lambda_i^k \right) + \beta \| \lambda_i^{k+1} - \lambda_i^k \|^2 \ge 0,\nonumber
\end{eqnarray}
where the last inequality is due to the definitions of $G_i$ and the assumptions that $\tau_i > \frac{\rho (A_i^T A_i)}{\beta}$ for inexact version. As a result, we have
\begin{eqnarray}
&& L(x_1^{k+1},\cdots,x_m^{k+1},\lambda) - L(x_1,\cdots,x_m,\lambda^{k+1})\nonumber\\
& \le & \sum\limits_{i=1}^m \left\{ \frac{1}{2} \left( \| x_i - x_i^k \|_{G_i}^2 - \| x_i - x_i^{k+1} \|_{G_i}^2 \right) + \frac{\beta}{2} \left( \| \lambda - \lambda_i^k \|^2 - \| \lambda - \lambda_i^{k+1} \|^2 \right)\right\}.\nonumber
\end{eqnarray}
Further summing from $k=0$ to $k=K-1$, we obtain
{\small\begin{eqnarray}
\sum\limits_{k=0}^{K-1} L(x_1^{k+1},\cdots,x_m^{k+1},\lambda) - \sum\limits_{k=0}^{K-1} L(x_1,\cdots,x_m,\lambda^{k+1})\le \sum\limits_{i=1}^m \left\{ \frac{1}{2} \| x_i - x_i^0 \|_{G_i}^2 + \frac{\beta}{2} \| \lambda - \lambda_i^0 \|^2 \right\}.\nonumber
\end{eqnarray}}
Letting $\tilde{x}_i^K := \frac{1}{K} \sum_{k=1}^K x_i^k,\quad \tilde{\lambda}^K := \frac{1}{K} \sum_{k=1}^K \lambda^k$, and using the convexity of the objective functions, we conclude
{\small$$
L(\tilde{x}_1^K,\cdots,\tilde{x}_m^K,\lambda) - L(x_1,\cdots,x_m,\tilde{\lambda}^K)\le \frac{\sum\limits_{i=1}^m \left\{ \| x_i \!-\! x_i^0 \|_{G_i}^2 + \beta \| \lambda \!-\! \lambda_i^0 \|^2 \right\}}{2K},\ \forall \!\left(\!\!\begin{array}{c} \mbox{\boldmath$x$}\\  \lambda \end{array}\!\!\right)\!\!\in {\mathcal{X}}\times{\mathbb{R}}^{\ell},
$$}
which is equivalent to
{small$$
\theta(\tilde{\mbox{\boldmath$x$}}^K) - \theta(\mbox{\boldmath$x$})
+ \!\left(\!\begin{array}{c}\tilde{\bfx}^K - \mbox{\boldmath$x$}\\ \tilde{\lambda}^K - \lambda \end{array}\!\right)\!^T F_2(\bfx,\lambda)\le \frac{\max\limits_{\bfx\in {\mathcal{X}},\lambda\in {\mathbb{R}}^{\ell} }\ \sum\limits_{i=1}^m \left\{ \| x_i \!-\! x_i^0 \|_{G_i}^2 + \beta \| \lambda \!-\! \lambda_i^0 \|^2 \right\}}{2K}.
$$}
The proof is complete. \hfill$\square$

\bibliography{Reference}
\bibliographystyle{plain}

\end{document}